\documentclass[twocolumn,amsthm]{autart}

\usepackage{empheq}

\usepackage{amsmath,amsfonts,amssymb,amscd,latexsym,enumerate,bbm}
\usepackage{array}
\theoremstyle{plain}
\newtheorem{theorem}{Theorem}[section]

\newtheorem{lemma}[theorem]{Lemma}
\newtheorem{corollary}[theorem]{Corollary}
\theoremstyle{definition}

\newtheorem{assumption}[theorem]{Assumption}
\newtheorem{remark}[theorem]{Remark}

\usepackage{arydshln}
\usepackage{tikz}
\usepackage{tkz-berge}

\usepackage{cite,url}
\usepackage{dsfont}

\usepackage{booktabs,siunitx}
\usepackage{mathtools}
\usepackage{dashbox}

\usepackage{enumerate}
\usepackage{bm}
\usetikzlibrary{calc,positioning,patterns,decorations.pathmorphing,decorations.markings,matrix,fit}
\tikzstyle{branch} = [circle,inner sep=0pt,minimum size=1mm,fill=black,draw=black]

\usepackage{tikz}
\usepackage{tkz-berge}
\usetikzlibrary{positioning}
\usetikzlibrary{arrows,%
	petri,%
	topaths}%
\usetikzlibrary{decorations.markings}         
\tikzstyle{vertex}=[circle, shading = ball, ball color = white!100!white, minimum size = 15pt, draw, inner sep=0pt]  

\DeclareMathOperator{\Tr}{Tr} 
\DeclareMathOperator{\vect}{vec} 
\DeclareMathOperator{\proj}{proj}
\DeclareMathOperator{\cl}{cl}
\DeclareMathOperator{\bd}{bd}
\DeclareMathOperator{\col}{col}

\DeclareMathOperator{\inter}{int}

\def\calb{{\mathcal B}}
\def\calc{{\mathcal C}}

\def\cali{{\mathcal I}}

\def\caln{{\mathcal N}}

\def\cals{{\mathcal S}}
\def\calt{{\mathcal T}}

\def\calx{{\mathcal X}}

\def\caly{{\mathcal Y}}

\newcommand{\agent}{\bm{\Sigma}}

\DeclareMathAlphabet{\mymathbb}{U}{bbold}{m}{n}

\newcommand{\bone}{\mathds{1}}
\newcommand{\bze}{\mymathbb{0}}

\newcommand{\R}{\mathbb{R}}
\newcommand{\N}{\mathbb{N}}

\newcommand{\ut}[2]{_#1{(#2)}}

\newcommand{\ud}{\, \mathrm{d}}

\DeclareMathOperator*{\argmin}{arg\,min}

\usepackage[prependcaption,colorinlistoftodos]{todonotes}

\usepackage[normalem]{ulem}

\usepackage{setspace}
\usepackage{url}
\usepackage[hidelinks]{hyperref}
\newcommand{\oprocendsymbol}{\hbox{$\bullet$}}
\newcommand{\oprocend}{\relax\ifmmode\else\unskip\hfill\fi\oprocendsymbol}

\allowdisplaybreaks

\newcommand{\BP}{\noindent{\bf Proof. }}
\newcommand{\EP}{\hspace*{\fill} $\blacksquare$\smallskip\noindent}

\journal{arXiv\date{}}

\begin{document}
	
	\begin{frontmatter}
		\title{ \LARGE 	\bf{Steering the aggregative behavior of noncooperative agents: a nudge framework}}
		\author[RUG]{Mehran~Shakarami}\ead{m.shakarami@rug.nl},$\,$
		\author[RUG]{Ashish Cherukuri}\ead{a.k.cherukuri@rug.nl},$\,$
		\author[RUG]{Nima Monshizadeh}\ead{n.monshizadeh@rug.nl}    
		\address[RUG]{Engineering and Technology Institute Groningen, University of Groningen, Nijenborgh 4, Groningen, 9747 AG , The Netherlands}
		
		\begin{abstract}
			
			This paper considers the problem of steering the aggregative behavior of a population of noncooperative price-taking agents towards a desired behavior. Different from conventional pricing schemes where the price is fully available for design,  we consider the scenario where a system regulator broadcasts a price prediction signal that can be different from the actual price  incurred by the agents. The resulting reliability issues are taken into account by including trust dynamics in our model, implying that the agents will not blindly follow the signal sent by the regulator, but rather follow it based on the history of its accuracy, i.e, its deviation from the actual price. We  present several nudge mechanisms to generate suitable price prediction signals that are able to  steer the aggregative behavior of the agents to stationary as well as temporal desired  aggregative behaviors. We provide analytical convergence guarantees for the resulting multi-components models. In particular, we prove that the proposed nudge mechanisms earn and maintain full trust of the agents, and the aggregative behavior converges to the desired one. The analytical results are complemented by a numerical case study of coordinated charging of plug-in electric vehicles.
		\end{abstract}
	\end{frontmatter}
	\endNoHyper
	\section{Introduction}
	Nudging is an approach in behavioral economics that is proposed to improve people's health and happiness by providing ``indirect suggestions'' termed as \emph{nudges}. A nudge, by definition,  is any characteristic of the choice structure that predictably changes people’s behavior without restricting any options or significantly affecting economic incentives\footnote{Nudge was originally defined as \textit{the minimalist intervention} in a given situation such that a desired outcome is achieved \cite{wilk1999mind}. However, the Nobel laureate Richard Thaler  presented another definition in \cite{thaler2008nudge} which is more popular and  is used here.}.  Therefore nudges are different from mandates as they are easy and cheap to avoid \cite{thaler2008nudge}. 
	Due to their aspects of preserving freedom of choice and being non-intrusive, nudge policies have become popular  over the last few years. The most notable example is the ``Behavioural Insights Team’’ (known as the ``Nudge Unit’’) that applies nudge theory in British government, and, for instance, its most recent report concerns energy consumption analysis and the impact of smart meters on customers’ energy consumption \cite{bit2020}. 
	Another example is  ``informational nudging'', defined as  sending manipulated, and possibly misleading, information about options to a decision maker for altering  its choices  \cite{guers2013informational}.  Informational nudging is studied recently in the context of 
	transportation systems \cite{cheng2016model} and boundedly rational decision makers \cite{cheng2018informational}.

	The problem of coordinating a population of noncooperative price-taking agents and altering their aggregative behavior appears in various applications 
	such as charging of plug-in electric vehicles in a coordinated way~\cite{ma2011decentralized},  
	residential energy consumption scheduling \cite{mohsenian2010optimal}, and congestion control in networks \cite{barrera2014dynamic}. 
	To address this problem, a common approach in the literature is treating the price as a design signal.
	If the system regulator has access to all information of the agents, a linear price with respect to the actions of the agents is sufficient to achieve a desired behavior \cite{alpcan2009nash}. In case such information is not available, which is often the case, dynamic pricing algorithms are posed as a solution to overcome this lack of knowledge; 
	see e.g.  \cite{alpcan2009nash,grammatico2016exponentially,ma2014distributed,barrera2014dynamic,grammatico2017dynamic}.
	The underlying assumption in dynamic pricing is that price is fully controllable, which  in turn facilitates the regulator's task in steering the behavior of the agents. 
	However, the actual price could depend on various elements such as fixed and variable production costs and daily market conditions; see e.g.  \cite{padhy2004unit} in the context of power systems. Here, instead, we allow the signal designed by the regulator to be different from the actual price dictating the costs incurred by the agents. Motivated by the advantages of nudging, 
	we propose a framework in which the regulator alters the aggregative behavior of price-taking agents, without directly designing the price and without fully knowing the cost/utility functions of the agents. In our setup, the regulator transmits a price prediction signal to all the agents. The agents choose their actions taking this prediction into account; however, they do not blindly follow it since they are aware that the prediction signal can differ from the actual price that they will incur. We model such behavior by associating a trust variable to each agent, which increases/decreases depending on the history of the accuracy of the communicated price prediction. 
	In other words, 
	here the agents cross-check the validity of the communicated information. This novel cross-checking step is a  key feature of our work, and distinguishes it from the existing informational nudging schemes \cite{guers2013informational,cheng2016model,cheng2018informational}.
	Moreover, the trust dynamics couple the price prediction dynamics to the actual price,  consequently the proposed nudge mechanisms do not simplify to conventional dynamic pricing schemes. 
	
	The presented framework is referred to as a nudge since  it does not directly affect economic incentives of the agents and respects their freedom of choice. Putting it differently, we use price information to indirectly suggest desired behaviors to the agents rather than enforcing mandates. 
		For the idea of nudging through price information in a different discipline, namely  agricultural economics, we refer the interested reader to \cite{belay2020nudging}.
	
	{\em Contributions:} We present a novel framework which is able to capture the multi-components model resulting from nudge mechanisms in conjunction with agents' actions and trust dynamics\footnote{Preliminary results of this work are presented in the conference article \cite{shakarami2020cdc}.
		Different to the conference article, this paper reports the proofs of Theorems \ref{thm:hard} and \ref{thm:soft}, studies convergence for stationary desired behaviors outside of the admissible set (Corollary \ref{cor:admi-hard}), presents a nudge mechanism for temporal desired behaviors (Section \ref{s:adaptive})  and establishes its convergence (Theorem \ref{t:new-adaptive} and Appendix \ref{apndx:uub}),  applies these results to coordinated charging of plug-in electric vehicles (Section \ref{sec:sim}), and studies existence of solutions for nonautonomous projected dynamical systems (Appendix \ref{app:pds-exist}).}.
	Within this framework, we first consider stationary desired behaviors and design two nudge mechanisms for the regulator, termed hard and soft nudge. We show that under these mechanisms, full trust of agents is gained in finite time and the aggregative behavior of the agents converges asymptotically to a desired set point. Afterwards, we extend the results to temporal desired behaviors and present an adaptive  nudge mechanism that can cope with the variations in the desired behavior. We analytically show that this mechanism obtains and maintains full trust of agents, and consequently the aggregative behavior converges to the time-dependent desired behavior. 
	Moreover, a byproduct of our analysis gives sufficient conditions for existence of Carath\'{e}odory solutions for nonautonomous
		projected dynamical systems.

	The structure of the paper is  as follows. Preliminaries are provided in Section~\ref{sec:math-pre}. The proposed framework is introduced in Section \ref{sec:prob_for}. Section~\ref{s:nudge} includes   the hard and soft  nudge mechanisms for stationary desired behaviors and their convergence analysis. The adaptive nudge mechanism for temporal desired behaviors is presented in Section \ref{s:adaptive}.  The case study is included in Section~\ref{sec:sim}, and finally,  conclusions are drawn in Section~\ref{sec:con}. Existence of solutions for nonautonomous projected dynamical systems is established in Appendix \ref{app:pds-exist} and stability analysis for the adaptive nudge is provided in Appendix \ref{apndx:uub}.

	\textbf{Notation.} We denote the set of natural, real, and nonnegative real numbers by $\N$,  $ \R $, and $ \R_{\geq0} $, respectively. The standard Euclidean norm is denoted by $\| \cdot \|$. The symbols  $ \bone_n $ and $ \bze_n $ respectively denote the vectors of all ones and zeros in $\R^n$. We denote the Kronecker product by $ \otimes $. The vectorization of a matrix $ M\in\R^{m\times n} $ is denoted by $ \vect(M) $. We denote the boundary, the interior, and the closure of a set $\calx \subseteq \R^n$ with $\bd(\calx)$, $ \inter(\calx) $, and $\cl(\calx)$, respectively. 
	Given 
	the vectors  $ x_1,\cdots,x_N \in\R^n$, we use the  notation $\col(x_i)=\big[x_1^\top ,\cdots,x_N^\top\big]^\top$. 
	We write $M\succ 0$ to indicate that $M=M^\top \in\R^{n\times n}$ is positive definite.  
	For a given vector $x\in \R^n$ and a positive semidefinite matrix $M$, we denote the weighted Euclidean norm of $ x $ by $ \|x\|_M:=\sqrt{x^\top M x} $. The Frobenius norm of a matrix $M\in\R^{m\times n}$ is denoted by $\|M\|_\text{F}:=\sqrt{\Tr(M^\top M)} $ where $\Tr(\,\cdot\,)$ is the trace operator. A closed ball with center $ x\in\R^n $ and radius $ r>0 $ is denoted by $ \bar{B}(x,r):=\{y\in\R^n\mid \|x-y\|\leq r\} $. 
	A  function $ F:\calx\to \R^m$ is locally Lipschitz on an open set $ \calx\subset \R^n $ if for any point $ x\in\calx $, there exist some positive scalar $ r$ and Lipschitz constant $ L $, both dependent on $ x $, such that $ \|F(y')-F(y)\|\leq L \|y'-y\| $ for all $y', y\in \bar{B}(x,r) $. The function $ F $ is Lipschitz on $ \calx $ if there exists a positive constant $L$ satisfying $\|F(y') - F(y)\| \le L \|y' - y\|$ for all $y', y \in \calx$.

	\section{Preliminaries}\label{sec:math-pre}
	This section provides basic notions on convex analysis and projected dynamical systems.
	
	\textit{Convex analysis:} Consider a nonempty, closed, convex set $\calx \subseteq \R^n$.  
	The map $ \proj_\calx:\R^n\to \calx $ denotes the Euclidean projection on to the set $ \calx $, i.e., $ \proj_\calx (z):=\argmin_{y\in\calx}\|y-z\| $. The normal cone to $ \calx $ at a given point $ x\in\calx $ is the set $ \caln_{\calx}(x):=\big\{y\in\R^n \mid y^\top(s-x)\leq 0,\forall s\in\calx\big\} $, and the tangent cone  is defined as the set $ \calt_\calx(x):=\cl\left(\cup_{y\in\calx}\cup_{\lambda>0} \lambda (y-x)\right) $. The projection of a vector $ z\in\R^n $ on to $ \calt_\calx(x) $ is denoted by $\Pi_{\calx}(x,z):= \proj_{\calt_\calx(x)}(z) $. Given any point $ x\in \calx $, it follows from Moreau's decomposition theorem \cite[Thm. 3.2.5]{hiriart1996convex} that any vector $ z\in \R^n $ can be written as $ z=\proj_{\caln_\calx(x)}(z)+\proj_{\calt_\calx(x)}(z) $. The  reader may refer to \cite[Fig. 5.3.1]{hiriart1996convex} for a geometrical representation of normal and tangent cones.
	
	\textit{Projected dynamical systems:} Given a nonempty closed set $\calx \subseteq \R^n$ and a continuous function $h:\R^n \times [0,\infty)\to\R^n$, the 
	non\-autonomous  projected dynamical system associated to them is 
	\begin{equation}\label{eq:pds}
		\dot x=\Pi_\calx  (x,h(x,t))\,.
	\end{equation} 
	The right-hand side of this system is discontinuous on the boundary of the set $\calx$. Following \cite[Def. 2.5]{nagurney2012projected}, we specify a notion of solution to the above projected dynamical system. 
	A map $x:[0,\infty)\to \calx$ is a Carath\'{e}odory solution of the projected dynamical system~\eqref{eq:pds} if it is absolutely continuous and satisfies  $\dot x(t)=\Pi_\calx (x(t),h(x(t),t))$ for almost all $ t \in [0,\infty)$.

	\section{Problem formulation and the model}\label{sec:prob_for}
	We consider a set of agents $ \cali := \{1, \dots , N\} $ that interact repeatedly with a central regulator. The agents are non\-cooperative, that is, each  agent $i$ is associated with a cost function $J_i$ that it wishes to minimize by choosing its action.  	In particular, the cost function of each agent $ i \in \cali$ is given by   $J_i\left(z_i,p\right)$, which determines the total
	cost of action $z_i\in\R^n $ given the price $ p\in\R^n $ and $n\in\N$. For simplicity, we assume that $J_i$ admits the following linear-quadratic form
	\begin{equation}\label{cost}
		\begin{aligned}
			J_i\left(z_i,p\right):=&\, \frac{1}{2}\left(z_i-c_i\right)^\top Q_i \left(z_i-c_i\right)+ z_i^\top p,
		\end{aligned}
	\end{equation}
	where 
	$ Q_i=Q_i^\top\in\R^{n\times n} $, $ Q_i\succ0 $, and $ c_i\in\R^n $. The cost function $ J_i $ consists of two terms, the local  penalty term $ \frac{1}{2}\left(z_i-c_i\right)^\top Q_i \left(z_i-c_i\right) $ and the cost of action $ {z_i}^\top p $. Note that $ c_i $ is the optimal action of the agent  when the price is zero.  The structure \eqref{cost} appears in applications where $z_i$ indicates the demand of a product that comes at price $p$, for instance coordinated charging of plug-in electric vehicles \cite{ma2011decentralized} and scheduling of residential energy consumption \cite{mohsenian2010optimal}.

	Before providing further details, we give an overview of our model. The regulator provides a prediction of the price for all the agents. This prediction is potentially different from the \emph{actual} price that determines the costs incurred by the agents. The agents use the price prediction to choose their actions with the aim of minimizing the cost they incur under the actual price. 
	The actual price is determined and revealed only after the actions are chosen.

	The regulator, on the other hand, aims at steering \textit{the aggregative behavior} of the agents to a desired point using the price prediction signal. We assume that the regulator does not know the cost functions of the agents.  A common approach of steering aggregate behavior, often referred to as dynamic pricing, is to use the price as a control signal to regulate the system of agents \cite{alpcan2009nash,grammatico2016exponentially,ma2014distributed}. In contrast, here the actual price signal is not available for design and the regulator needs to rely on the price prediction signal to manipulate the agents' behavior. Our motivation stems from the fact that, in reality, the actual price may not be prescribed a priori  as a dynamic function of demands/actions.
	
	{The discrepancy between the price prediction and the actual price} readily brings the issue of {\em trust} or \textit{reliability}. Namely, the central regulator needs to earn and maintain the trust of the agents in order to influence their decisions. 
	We take this into account by considering that the agents associate a level of trust/reliability to the regulator's prediction based on the \textit{history of its accuracy}.
	
	In the sequel, we aim to carefully model the above described features and design update schemes, termed \textit{nudge mechanisms}, that enable the regulator to steer the aggregative behavior of the agents to a desired reference. 
	We first look at the problem from the agents' side and put forward a model where agents use available information to decide on their actions. 
	The regulator's side will be dealt with in Section \ref{s:nudge}, where  nudge mechanisms are proposed. 
	\subsection{Agents' actions and trust dynamics}
	In choosing their actions at time $t \in [0,\infty)$, the agents have access to a price prediction $ \hat{p}{(t)}\in \R^n$ sent out by the regulator. Note that this value is common for all agents. In addition, we assume that each agent $i \in \cali$ has a local perception of the price, denoted by $ \hat{\lambda}_i\in \R^n$, that the agent would have used in the absence of the prediction $ \hat{p}{(t)} $. 
	
\begin{figure}
		\centering
		\includegraphics[width=0.6\columnwidth]{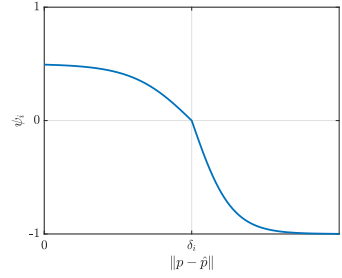}
		\caption{An example of the function $ \psi_i $ satisfying Assumption \ref{asm:psi}.}
		\label{fig:psi}
	\end{figure}

	As mentioned before, different from conventional dynamic pricing,  the distinction between the actual price and its prediction brings the issue of \textit{reliability}, and we incorporate this in our model by associating a level of trust/reliability to the regulator's prediction based on the history of its accuracy. In particular, let $ \gamma\ut{i}{t}\in [0,1] $ be the trust variable of agent $ i $ associated with the price prediction $ \hat{p}{(t)} $. Note that $\gamma_i(t) =1$ and $\gamma_i(t) = 0$ stand for full and no trust, respectively. 
	Given the amount of trust, predicted price, and the local perception, agent $i$ 
	adopts a \textit{trust-adapted} 	price perception\footnote{\label{fn:trust}The trust-adapted protocol \eqref{local_price} can be replaced by a more general form 
				$\lambda_i(t) = \omega_i( \hat p (t), \gamma_{i}(t), t)$ where $(\hat{p},\gamma_i)\mapsto \omega_i( \hat{p},\gamma_i, t)$ is  Lipschitz, $ t\mapsto  \omega_i( \hat{p},\gamma_i, t)$ is uniformly continuous, and   $\omega_i(\hat p,1,t)=\hat p$ for all $ t\in [0,\infty) $. The explicit dependence of $\omega_i$ on $t$ also allows to accommodate a time-varying local price perception $t\mapsto\hat \lambda_i(t)$. However, we opt for the form \eqref{local_price} in order to provide a more explicit analysis and to highlight better the underlying intuition.}
	\begin{equation}\label{local_price}
		\lambda\ut{i}{t}:=\gamma\ut{i}{t}\hat{p}{(t)}+(1-\gamma\ut{i}{t})\hat\lambda_i\,.
	\end{equation}
	If $\gamma_i(t)$ is close to $ 1 $, the agent disregards its own perception of the price and follows the price prediction communicated by the regulator. Conversely, as $\gamma_i(t)$ approaches $ 0 $, the agent loses trust in 
	the price prediction $\hat{p}(t)$ and follows its own price perception $\hat \lambda_i$ when deciding on its optimal action.  	The agent $i$ uses this trust-adapted 
	price 	perception to determine its optimal action as follows:
	\begin{equation*}
		x\ut{i}{t}:=\argmin_{z\in\R^n}J_i\left(z,\lambda\ut{i}{t}\right).
	\end{equation*}
	By using \eqref{cost} and \eqref{local_price}, the 
	explicit expression of the optimal action of agents is given by
	\begin{equation}\label{opt_action}
		x\ut{i}{t}=c_i-Q_i^{-1} \Bigl(\gamma\ut{i}{t}\hat{p}{(t)}+\left(1-\gamma\ut{i}{t}\right)\hat\lambda_i\Bigr).
	\end{equation}
	The actual price $t\mapsto p(t)$ is available to the agents once they have taken their actions. 	If the discrepancy between the predicted and actual price is large, then agents lose their trust in the predictions. 	We capture the changes of trust based on these positive or negative experiences by providing a trust update rule. 
	In particular, we consider the following trust dynamics:
	\begin{equation}\label{trust_core}
		\dot \gamma_i(t)=\eta_i \psi_i(\|p(t)-\hat{p}{(t)}\|), 
	\end{equation}
	where $\eta_i>0 $ and $ \psi_i:\R_{\geq 0}\to [-1,1]$ determines whether the agent loses or gains trust  in the price prediction. We assume that $ \psi_i(\, \cdot \,) $ satisfies the following assumption, and an example of 
	this function is depicted in Fig.~\ref{fig:psi}.

	\begin{assumption}\label{asm:psi}
		The  function $ \psi_i:\R_{\geq 0}\to [-1,1]$ is   locally Lipschitz     and strictly decreasing. In addition, we have 
		$ \psi_i(0)>0 $ and $ \psi_i(\delta_i)=0 $ for some $ \delta_i>0 $. \oprocend
	\end{assumption}
	
	The scalar $ \delta_i $ quantifies the \textit{tolerance} of agent $ i $ towards the prediction error. That is, if the error between the actual and the predicted price $\|p(t)-\hat{p}{(t)}\|$ is greater than $ \delta_i $, agent $ i $ begins losing trust in the prediction with the rate $\eta_i$. Conversely, trust increases as long as the error 
	is within the tolerance $\delta_i$. The rationale behind this dynamics is that, excluding the extreme cases of unconditional trust or distrust, trust  can be gained or lost after several positive or negative experiences \cite{jonker1999formal}.

	Note that trust variables are defined in the interval between $0$ and $1$.
	To respect this, we slightly revise \eqref{trust_core} by adding projection operators to it, namely:
	\begin{equation}\label{trust_proj}
		\dot \gamma_i(t)=\Pi_{[0,1]}\left(\gamma_i(t),\eta_i \psi_i(\|p(t)-\hat{p}{(t)}\|)\right).
	\end{equation}
	We note that the essence of the trust update rule remains the same as \eqref{trust_core}. The projection operators become active only if the bounds $\gamma_i=0$ or $\gamma_i=1$ are hit. In particular, if $\gamma_i(t_1)=1$ at some time $t=t_1$ and $\psi_i(\|p(t_1)-\hat{p}{(t_1)}\|)$ is positive (thus suggesting an increase in $\gamma_i$), the projection becomes active, and sets $\dot \gamma_i(t_1)$ to $0$, thus prohibiting the trust variable to exceed its maximum value $1$. An analogous scenario occurs for the case $\gamma_i(t_1)=0$.

	For simplicity of presentation, we rewrite the model of agent $ i $, consisted from \eqref{opt_action} and \eqref{trust_proj}, as follows:
	\begin{subequations}\label{agent}
		\begin{empheq}[left={\agent_i:\empheqlbrace}]{align}
			\dot\gamma_i(t)&=\Pi_{[0,1]}(\gamma_i(t),\eta_i \psi_i(\|p(t)-\hat{p}{(t)}\|)),\label{trust_dyn}\\
			x\ut{i}{t}&=\pi_i({\hat{p}}(t), \gamma_i(t))\label{opt_action_cont},
		\end{empheq}
	\end{subequations}
	where 
	\begin{equation}\label{pi_def}
		\pi_i({\hat{p}}, \gamma_i):=c_i-Q_i^{-1} \Bigl(\gamma_i\,\hat{p}+\left(1-\gamma_i\right)\hat\lambda_i\Bigr).
	\end{equation}
	Note that the actual price $ p $ and the price prediction $ \hat{p} $ are the inputs of the  model, and the action vector $ x_i $ is the output. 
	Having introduced the model of the agents, we next discuss the desired aggregative behavior.

	\subsection{Desired aggregative behavior}
	The goal of the system regulator is to coordinate the  agents such that they cumulatively behave in a desired fashion.  Here, we are interested in regulating  $\sum_{i\in\cali}x_i(t)  $, which we refer to as the \textit{aggregative behavior}. Such quantity often reflects total production or total demand depending on the application at hand. 
	More precisely,  the regulator aims to achieve 
	\begin{equation}\label{e:reg-goal}
		\lim\limits_{t\to \infty}\sum_{i\in\cali}x\ut{i}{t}=x^*,
	\end{equation}
	for some desired setpoint $ x^*\in\R^n $.\footnote{In Section \ref{s:adaptive}, we allow $x^*$ to be a time-varying reference signal.} To this end, we propose suitable \textit{nudge} mechanisms that can be implemented by the regulator.  A mechanism is  a nudge if it influences the behavior of a group of individuals through providing indirect suggestions. We use this concept and propose  mechanisms in which the  regulator manipulates the price prediction $ \hat{p}{(t)} $ to achieve its goal, namely \eqref{e:reg-goal}.

	Recall that the actual price is considered here as an exogenous signal. In particular, we assume that it admits 
	\begin{equation*}
		p(t)=p_0+\Delta p(t)\,,\quad \forall\,  t \in [0,\infty)\,,
	\end{equation*}
	where $ p_0 $ is a constant base price, known to the regulator, and $ \|\Delta p(t)\| \ll \|p_0\|$ accounts for price fluctuations. We assume that the following condition holds throughout the paper:

	\begin{assumption}\label{std_asm}
		The actual price function  $p : [0,\infty) \rightarrow \R^n$ is continuous, 
		and its fluctuations satisfies 
		$  \|\Delta p(t)\|< \min_{ i\in \cali  }\delta_i$ for all $ t\in [0,\infty) $. \oprocend
	\end{assumption}

	\begin{remark}\label{rem:delta-p}
		Note that in the absence of the objective \eqref{e:reg-goal}, the best the regulator can do is to provide the agents with the true value of $p_0$. In that case, the price prediction error amounts to $\|\Delta p(t)\|$. Therefore, the inequality constraint in Assumption \ref{std_asm} simply means that the prediction error in such a \textit{manipulation-free} case is within the tolerances of all agents. In other words, the price fluctuations, per se, should not lead to a loss in trust. \oprocend
	\end{remark}

	The fact that the agents do not blindly follow the price prediction $ \hat{p}{(t)} $ implies that not any arbitrary aggregative behavior $x^*$ is achievable.  
	Next, we identify a set of aggregative behaviors to which the agents can be driven by applying our nudge mechanisms.

	Let Assumption~\ref{std_asm} hold, and choose $ \bar \delta\in \R$ such that
	\begin{equation}\label{delta_bar}
		0<\bar \delta< \min_{i\in\cali}\delta_i- \|\Delta p(t)\|,\quad \forall\, t\in [0,\infty).
	\end{equation}
	We leverage the idea that if Assumption \ref{asm:psi} holds and $ \hat{p}{(t)} $ belongs to the closed ball
	\begin{equation}\label{ball}
		\calb:=\bar{B}(p_0,\bar\delta)=\left\{\hat{p}\in\R^n\mid \|\hat{p}-p_0\|\leq \bar\delta\right\},
	\end{equation}
	then $\psi_i(\,\cdot\,)$ takes positive values and $\gamma\ut{i}{t} $ increases for all $ i\in\cali$ following \eqref{trust_dyn}. As a result, the regulator can gain agents' trust in the price prediction by constraining $ \hat{p}{(t)} $ to the ball $  \calb$. 	
	Bearing this and the action of agents in \eqref{opt_action_cont}  in mind, we define the set of admissible $x^*$ as: 
	\begin{equation}\label{xs_def}
		\calx^*:=\Big\{x\in\R^n\mid x=\sum_{i\in\cali}\left(c_i-Q_i^{-1}\hat{p}\right), \, \hat{p}\in \calb\Big\}.
	\end{equation}
	From~\eqref{ball}, the set $\calx^*$ can be explicitly written as
	\begin{equation}\label{xs_def_explicit}
		\calx^*=\Big\{x\in\R^n\mid (x -x_0)^\top (\sum_{i\in\cali}Q_i^{-1})^{-2} (x-x_0) \le \bar \delta^{\,2} \Big\},
	\end{equation}
	where $x_0:=\sum_{i\in\cali}\left(c_i-Q_i^{-1}{p}_0\right)$. Thus, the regulator can alter the aggregative behavior inside a compact set around $x_0$. 
	Putting it differently, $\calx^*$ characterizes the set of aggregative behaviors that are potentially achievable while monotonically increasing the trust variables. Note from~\eqref{delta_bar} and  \eqref{xs_def_explicit} that the bigger the agents' tolerances $\delta_i$'s are, the larger can be $\bar \delta$ and thus, the  set $\calx^*$.

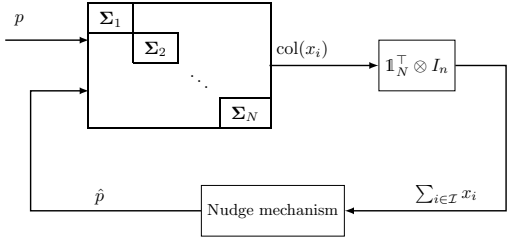
\begin{figure}
		\centering
		\resizebox{0.8\columnwidth}{!}
		{
			\begin{tikzpicture}
				\node (agents)  {
					\begin{tabular}{|llll|}
						\hline
						\multicolumn{1}{|l|}{$\agent_1$} &                       &                       &                         \\ \cline{1-2}
						\multicolumn{1}{|l|}{} & \multicolumn{1}{l|}{$\agent_2$} &                       &                        \\ \cline{2-2}
						&                       &               $ \ddots $          &                     \\ \cline{4-4}
						&                       & \multicolumn{1}{l|}{} & \multicolumn{1}{l|}{$\agent_N$} \\ \hline
				\end{tabular}};
				\node (Kron) [draw,rectangle,minimum width=1.5cm, minimum height=1cm, right=2cm of agents] {$ \bone_N^\top\otimes I_n$};
				\node (nudge) [draw,rectangle,minimum width=1.5cm, minimum height=1cm,  below right=1cm and -1.5cm of agents] {\footnotesize Nudge mechanism};

				\draw[-latex, thick]        ($(agents.west)+(-1.5,0.5)$)-- ++(1.65,0); 
				
				\draw[-latex, thick]       (agents.east)+(-0.15,0) -- (Kron.west);
				\draw[-latex, thick]        (Kron.east)-- ++(1,0)|- (nudge.east);
				\draw[-latex, thick]       (nudge.west) -| ($ (agents.west)+(-1,-0.5) $) -- ++(1.15,0);
				
				\draw[-latex,thick] node[right=0.5cm of agents,above] {$\col(x_i)$} node[right=2cm of nudge,above] {$\sum_{i\in\cali} x_i$} node[left=2cm of nudge,above] {$\hat{p}$} node[above left=-0.7 and 1cm of agents] {$p$};
		\end{tikzpicture}}
		\caption{Block diagram representation of agents  interconnected with a nudge mechanism.}
		\label{fig:block}
	\end{figure}

	For any $ x^*\in \calx^* $, there exists a unique $ p^*\in \calb $ such that
	\begin{equation}\label{x_star}
		x^*=\sum_{i\in\cali}\left(c_i-Q_i^{-1}{p}^*\right)\,,
	\end{equation}
	or equivalently
	\begin{equation}\label{p_star}
		p^*=\big(\sum_{i\in\cali}  Q_i^{-1} \big) ^{-1} \big(-x^*+ \sum_{i\in\cali} c_i\big)\,.
	\end{equation}
	The vector $p^*$ is an important quantity. If the agents fully trust the price prediction and the regulator communicates $p^*$ as the prediction, then the aggregative behavior of the agents will be $x^*$.   
	However,  the regulator cannot directly compute $p^* $ since it does not know the exact
	parameters defining individual cost functions. Moreover, trust can only be gained over time. To address these issues, suitable nudge mechanisms are designed in the next section. Each of those mechanisms can be interconnected with the agents' dynamics, as demonstrated in Fig.~\ref{fig:block}, in order to drive
	the price prediction $\hat{p}(t)$ to $p^*$, and consequently $x(t)$ to $x^*$.
	\label{edits-for-delta}
	The key parameter used in the proposed mechanisms is $ \bar{\delta} $  satisfying \eqref{delta_bar}. The  precise values of the tolerances of the agents $ \delta_i $'s are unknown to the regulator, and the price fluctuations $ \Delta p(t) $ are not available a priori. Thus the regulator typically needs to rely on lower estimate of $\min_{i\in\cali}\delta_i- \|\Delta p(\cdot)\|$ to select $ \bar{\delta} $. The less the regulator knows about the right-hand side of \eqref{delta_bar}, the more conservative the value of $\bar \delta$ has to be chosen, which in turn results in a  smaller ball $\mathcal{B}$ as well as a smaller set of admissible desired behaviors $ \calx^* $. Learning a feasible $\bar \delta$ from experiments is an interesting research question for future research.

	\section{Nudge mechanisms for stationary desired behaviors}\label{s:nudge}
	In this section, we design two nudge mechanisms, referred to as \emph{hard} and \emph{soft}, that provide suitable price prediction signals. 
	
	\subsection{Hard nudge mechanism}
	The first nudge mechanism that we propose is the following projected-integral control law
	\begin{equation}\label{hard_nudge}
		\dot{\hat{p}}(t)=\Pi_{\calb}\Big(\hat{p}(t),\sum_{i\in\cali}x_i(t)-x^*\Big),
	\end{equation}
	where $ \calb $ is defined as \eqref{ball}  
	and $x^*$ is the desired aggregative behavior. 
	We note that from~\cite[Lem. 2.1]{nagurney2012projected}, the projection operator on the right-hand side can be explicitly expressed using the definition of $ \calb $. In particular, let $ e(t):=\sum_{i\in\cali}x_i(t)-x^* $, then we obtain:
	\begin{equation*}
		\Pi_{\calb}\left(\hat{p}(t),e(t)\right)=\begin{cases}
			e(t), &  \text{if } \hat{p}(t)\in \inter(\calb),\\
			e(t)-\frac{\alpha(t)(\hat{p}(t)-p_0)}{\|\hat{p}(t)-p_0\|^2},& \text{if }\hat{p}(t)\in \bd(\calb),
		\end{cases}
	\end{equation*}
	where $ \alpha(t):=\max\{0,e(t)^\top (\hat{p}(t)-p_0)\} $.
	The intuition behind the nudge mechanism in \eqref{hard_nudge} is as follows: this mechanism provides a suitable  integral action that updates the price prediction such that the error between the desired behavior and the current aggregative behavior  diminishes. To gain and maintain the trust of the agents, the price prediction is constrained to the ball $\calb$ for all time, and thus we refer to \eqref{hard_nudge} as {\em hard} nudge.

	The overall system, as shown in Fig.~\ref{fig:block}, is obtained by interconnecting \eqref{hard_nudge} with  agents  
	\eqref{agent}, and the theorem below addresses its convergence.
	
	\begin{theorem}\label{thm:hard}
		Consider the closed-loop system formed by agents' model \eqref{agent} and the hard nudge mechanism  \eqref{hard_nudge} with   $ x^*\in \calx^* $.
		Then, for any initial condition $\left(\hat{p}{(0)},\col (\gamma_i{(0)})\right)\in {\calb}\times [0,1]^{N} $, there exists a  Carath\'{e}odory solution  
		$t \mapsto \left(\hat{p}{(t)},\col (\gamma_i{(t)})\right)$ of the closed-loop system over the domain $ [0,\infty) $. Moreover, any solution $ \left(\hat{p}{(t)},\col (\gamma_i{(t)})\right) $ 
		converges to $ (p^*,\bone_N) $ with $ p^* $ given by \eqref{p_star}. Consequently, $\sum_{i\in \cali} x_i(t)$  converges to $x^*$ as desired.
	\end{theorem}

	\BP
	The proof is divided into two parts. Since the vector field of the overall system is discontinuous, we  show existence of Carath\'{e}odory solutions of the system in the first part. The second part is devoted to convergence analysis.

	\textit{Existence of solutions: }
	Let $ \xi:=({\hat{p}},\col( \gamma_i)) $ and $ \Omega:=\calb\times [0,1]^{N} $. Then, by substituting the expression of $ x_i $ from  \eqref{opt_action_cont} into \eqref{hard_nudge}, we obtain the nonautonomous projected dynamical system that represents the closed-loop system \eqref{agent} and \eqref{hard_nudge} as 
	$
	\dot\xi=\Pi_{\Omega}(\xi,h(\xi,t))
	$,
	where
	\begin{equation*}
		h(\xi,t):=\begin{bmatrix}
			\sum\limits_{i\in\cali}\pi_i({\hat{p}}, \gamma_i)-x^*\\
			\col\left(\eta_i \psi_i(\|p(t)-\hat{p}\|)\right)
		\end{bmatrix}.
	\end{equation*}
	Note that the map $ (\hat{p},t) \mapsto \psi_i(\|p(t)-\hat{p}\|) $ is measurable\footnote{\label{foot:measur}A function $ f:E\to \R $ is measurable if its domain $ E $ is measurable, and the set $ \{x\in E \mid f(x)>\alpha \} $ is measurable for all $ \alpha \in \R $. For an in depth overview of measurable functions see \cite[Ch. 3]{royden1988real}.} in $t$ and locally Lipschitz  in $\hat{p}$. The former follows from Assumptions~\ref{asm:psi} and \ref{std_asm} and the fact that every continuous function is measurable \cite[Prop. 3.3]{royden1988real}.  The latter is a consequence of Assumption~\ref{asm:psi} and the fact that the norm operator is Lipschitz. 
	Consequently,  the function $ (\xi,t) \mapsto h(\xi,t) $ is locally Lipschitz in $\xi$ and measurable in $t$, and using the compactness of the set $\Omega$, existence of solutions for any initial condition $\left(\hat{p}{(0)},\col (\gamma_i{(0)})\right)\in {\calb}\times [0,1]^{N} $ is guaranteed by Lemma~\ref{lem:exis-pds}.

	\textit{Convergence analysis:} Our proof proceeds by showing that for any solution of the system, there exists a finite time by which full trust of agents is achieved and maintained. Subsequently, with full trust, we show that $\hat p(t)$ converges to $p^*$.
	
	Note from~\eqref{hard_nudge} that $ \hat{p}{(t)} \in \calb $ for all $t \geq 0$. Using this fact along with Assumption~\ref{asm:psi}, we obtain $ \psi_i(\|p(t)-\hat{p}(t)\|)>0 $ for all $ i\in\cali $ and $ t\geq 0 $. Consequently, along any solution, the trust variable of agent $i$ at any time $t$ is given by 
	\begin{equation}\label{trust_time}
		{\gamma}\ut{i}{t}=\min\Big\{1,\gamma\ut{i}{0}+\eta_i \!\int_{0}^{t}\!\!\psi_i(\|p(\tau)-\hat{p}(\tau)\|)\ud \tau\Big\}.
	\end{equation} 
	Bearing in mind that $\hat p(t)$ belongs to the ball $\calb$ given by \eqref{ball}, we have $\|p(t)-\hat{p}(t)\| \leq \rho <\min_{i\in\cali}\delta_i$ for some $\rho>0$. Hence, by Assumption~\ref{asm:psi}, we obtain that $\psi_i(\|p(t)-\hat{p}(t)\|)\geq \psi_i(\rho)>0$ for all time.  
	Let $ T^i:=(1-\gamma\ut{i}{0})/(\eta_i\psi_i(\rho)) $. Then, from~\eqref{trust_time}, we deduce that  $ {\gamma}\ut{i}{t}=1 $ for  all $ t\geq T^i $. Setting $ T:=\max_{i\in\cali}T^i $, we conclude that   $ \col(\gamma_i{(t)})=\bone_N $ for all $ t\in [T,\infty) $. 
	As a consequence, in the time interval $[T,\infty)$, the price prediction dynamics~\eqref{hard_nudge} reduces to 
	\begin{equation}\label{p_hat_dot}
		\dot{\hat{p}}=\Pi_{\calb}(\hat{p},f\left(\hat{p}\right)),
	\end{equation}
	where 
	\begin{equation}\label{f_map}
		\begin{aligned}
			f\left(\hat{p}\right):=\sum_{i\in\cali}c_i-\sum_{i\in\cali}Q_i^{-1}\hat{p}-x^*.
		\end{aligned}
	\end{equation}
	We next analyze the asymptotic properties of~\eqref{p_hat_dot} and show that its solutions  converge asymptotically to $p^*$.  
	Consider the Lyapunov candidate $V(\hat{p}):=\frac{1}{2}\|\hat{p}-p^*\|^2$. Since solutions of~\eqref{p_hat_dot} are absolutely continuous and $V$ is continuously differentiable, the time-derivative of the evolution of $V$ along any solution of~\eqref{p_hat_dot} is equal to the inner product of the gradient of $V$ and the right-hand side of~\eqref{p_hat_dot}.  
	This inner product is computed as 
	 \begin{align*}
		         \nabla V(\hat{p})^\top \Pi_{\calb}\left(\hat{p},f(\hat{p})\right) =&(\hat{p}-p^*)^\top f(\hat{p})\\
		        &-(\hat{p}-p^*)^\top\proj_{\caln_{\calb}(\hat{p})}\left(f\left(\hat{p}\right)\right),
	 \end{align*} 
	where we used Moreau’s decomposition theorem (cf. Section~\ref{sec:math-pre}) to obtain the above equality and $ \caln_{\calb}(\hat{p}) $ is the normal cone of $ \calb $ at $ \hat p $. Note that $ -(\hat{p}-p^*)^\top \proj_{\caln_{\calb}(\hat{p})}\left(f\left(\hat{p}\right)\right)\leq 0$ since $  \hat{p},p^*\in  \calb$, and we find that 
	$
	\nabla V(\hat{p})^\top \Pi_{\calb}\left(\hat{p},f(\hat{p})\right)\leq  (\hat{p}-p^*)^\top f(\hat{p})\,.
	$
	We use \eqref{f_map} and the expression of $ x^* $ in \eqref{x_star} to obtain
	\begin{equation*}
		\nabla V(\hat{p})^\top \Pi_{\calb}\left(\hat{p},f(\hat{p})\right)\leq - \|\hat{p}-p^*\|_{\sum_{i\in\cali}Q_i^{-1}}^2\,.
	\end{equation*}  
	This implies that $V$ decreases monotonically along every solution of~\eqref{p_hat_dot}. Consequently, $\hat p$ converges to $p^*$, and the aggregate behavior $\sum_{i\in \cali} x_i$ converges to $x^*$.
	\EP

	As shown in Theorem \ref{thm:hard}, the hard nudge mechanism \eqref{hard_nudge} successfully steers the agents to the desired aggregative behavior, for any $x^*\in \calx^*$. 
	\label{x-star-not-in-set}
	An implicit requirement is that  the regulator has partial knowledge on expected desired aggregative behaviors, i.e., a subset of $ \calx^* $, to pick a feasible $ x^* $. In case this information is not available and
	$x^*\notin \calx^*$, convergence of the aggregative behavior is still guaranteed, but to a different
	point, namely to $x^\prime \in \mathcal{X}^*$  that is the closest point to $x^*$ in a suitable norm.
	This is formally stated in the following corollary.
	
	\begin{corollary}\label{cor:admi-hard}
		Consider the closed-loop system formed by agents' model \eqref{agent} and the hard nudge mechanism  \eqref{hard_nudge} with   $ x^*\notin \calx^* $.
		Then, for any initial condition $\left(\hat{p}{(0)},\col (\gamma_i{(0)})\right)\in {\calb}\times [0,1]^{N} $, there exists a  Carath\'{e}odory solution 
		$t \mapsto \left(\hat{p}{(t)},\col (\gamma_i{(t)})\right)$ of the closed-loop system over the domain $ [0,\infty) $. Moreover,  
		$\sum_{i\in \cali} x_i(t)$ converges to $x'\neq x^*$ given by
		\begin{equation*}
			x'=\argmin_{y\in \calx^*}\frac{1}{2}\|x^*-y\|_{\left(\sum_{i\in\cali}Q_i^{-1}\right)^{-1}}^2\,.
		\end{equation*}
	\end{corollary}
	\BP
	Based on the proof of Theorem \ref{thm:hard}, the closed-loop system admits a Carath\'{e}odory solution for all $ x^*\in\R^n $, and thus existence of a solution $t \mapsto \left(\hat{p}{(t)},\col (\gamma_i{(t)})\right)$ is guaranteed for all $ t\in[0,\infty) $.  Next, we consider $ x' $ and characterize its corresponding  price prediction, namely $ p' $. We prove convergence of  $ (\hat{p},\col(\gamma_{i})) $ to $ (p',\bone_N) $ afterwards. Subsequently, convergence of $\sum_{i\in \cali} x_i$  to $x'$ follows from the definition of $ p' $.
	
	The point $ x' $ exists and is unique following  Weierstrass’ Theorem \cite[Prop. A.8]{bertsekas99nonlinear} and \cite[Prop. 2.1.1]{bertsekas99nonlinear}, respectively.  
	It also follows from \cite[Prop. 2.1.2]{bertsekas99nonlinear} that $ x'\in \calx^* $ satisfies 
	\begin{equation*}
		\left(x'-x^*\right)^\top \Big(\sum_{i\in\cali}Q_i^{-1}\Big)^{-1}\left(y-x' \right)\geq 0,\quad \forall\, y\in \calx^*.
	\end{equation*}
	Let $ p':=(\sum_{i\in\cali}  Q_i^{-1} ) ^{-1} (-x'+ \sum_{i\in\cali} c_i) $, then we have $ p'\in\calb $. Moreover, for all $y\in \calx^*$, we have
	\begin{equation*}
		\!\Bigl(\sum_{i\in\cali}(c_i-Q_i^{-1} p')-x^*\Bigr)\!{\Big.^\top}  \! \Bigl(\big(\sum_{i\in\cali}Q_i^{-1}\big)\!{\big.^{-1}}(y-\!\sum_{i\in\cali} c_i)+p' \Bigr)\!\geq 0.
	\end{equation*}
	Recalling the definition of $ \calx^* $ given by \eqref{xs_def}, we see that for any $ y\in\calx^*$, there exists some $ s\in\calb $ such that the relation $ y=\sum_{i\in\cali}\left(c_i-Q_i^{-1}s\right) $ holds. Therefore, the above inequality can be rewritten as  
	\begin{equation}\label{adms_p_ineq}
		\Big(\sum_{i\in\cali}\big(c_i-Q_i^{-1} p'\big)-x^*\Big)^\top\left( p'-s \right)\geq 0,\quad \forall\, s\in \calb.
	\end{equation}	
	Note from~\eqref{hard_nudge} that $ \hat{p}{(t)} \in \calb $ for all $t \geq 0$. Following the steps of the proof of Theorem \ref{thm:hard},  there exists some finite time $ T\geq 0 $ such that  $ \col(\gamma_i{(t)})=\bone_N $ and the hard nudge mechanism reduces to  \eqref{p_hat_dot} for all $ t\geq T $. Considering again the Lyapunov candidate $V(\hat{p}):=\frac{1}{2}\|\hat{p}-p'\|^2$, its derivation along \eqref{p_hat_dot} yields
	\begin{equation*}
		\nabla V(\hat{p})^\top \Pi_{\calb}\left(\hat{p},f(\hat{p})\right)\leq  (\hat{p}-p')^\top f(\hat{p}).
	\end{equation*}
	Now we add the left-hand side of  \eqref{adms_p_ineq} evaluated at $ s=\hat{p} $ to the right-hand side of   the foregoing inequality to get
	\begin{equation*} 
		\begin{split}
			\nabla V(\hat{p})^\top &\Pi_{\calb}\left(\hat{p},f(\hat{p})\right)\\&\leq  (\hat{p}-p')^\top \Big(f(\hat{p})-\sum_{i\in\cali}\big(c_i-Q_i^{-1} p'\big)+x^*\Big)\\
			&= - \|\hat{p}-p'\|_{\sum_{i\in\cali}Q_i^{-1}}^2,
		\end{split}
	\end{equation*}
	where the  equality follows from the definition of  $ f $ given by \eqref{f_map}. We conclude that $ V $ decreases monotonically along every solution of \eqref{p_hat_dot} and $ \hat{p} $ converges to $ p' $. 
	\EP
		\begin{remark}\label{rem:no-full-trust}
		If  Assumption \ref{std_asm} is not satisfied,  one may still be able to provide convergence guarantees under suitable conditions.
		In particular, let $\cals$ denote the collection of agents that violate Assumption \ref{std_asm} for all time, i.e, $
		\cals:=\{j\in \cali \mid  \|\Delta p(t)\|\geq \delta_j,\, \forall t\in [0,\infty)\}.
		$
		The remaining agents satisfy the assumption, namely $ \|\Delta p(t)\|<\min_{i\in\cali\setminus\cals}\delta_i $ for all time.
		We can then show that under the hard nudge \eqref{hard_nudge} with	$ \bar{\delta}\in\R $  satisfying the revised inequality
		\begin{equation*}
		0<\bar \delta< \min \big\{\min_{i\in\cali \setminus \cals}\delta_i- \|\Delta p(t)\|, \ \|\Delta p(t)\|-\max_{j\in \cals}\delta_j\big\}, 
		\end{equation*}
		for all time, the aggregative behavior of the agents in $\cals$ converges to $ \bar{x}:=\sum_{j\in \cals}(c_j-Q_j^{-1}\hat{\lambda}_j)$, whereas the aggregative behavior of the agents in $\mathcal{I}\setminus \cals$ converges to 
			\begin{equation*}
			x'=\argmin_{y\in \caly}\frac{1}{2}\|x^*-\bar{x}-y\|_{(\sum_{i\in\cali\setminus \cals}Q_i^{-1})^{-1}}^2\,,
			\end{equation*}
			where
			$
			\caly:=\{x\in\R^n\mid x=\sum_{i\in\cali\setminus \cals}\left(c_i-Q_i^{-1}\hat{p}\right), \, \hat{p}\in \calb\}
			$.  
			The set $\caly$ is similar to $\mathcal{X}^*$ in \eqref{xs_def} but with the set of agents restricted to $\mathcal{I}\setminus \cals$. In case $x^*-\bar{x} \in \caly$, we have $x'=x^*- \bar x$ which implies that the aggregative behavior of all agents converges to $x'+\bar x=x^*$. 
			The details of the analysis are omitted due to lack of space.
			\oprocend	
	\end{remark}

	\subsection{Soft nudge mechanism}\label{subsec:soft}
	While using  the nudge mechanism in \eqref{hard_nudge} is effective for driving the aggregative behavior of the agents to a desired point, convergence is guaranteed only if the price prediction is initialized in the ball $\calb$. We now present an alternative nudge mechanism under which convergence is guaranteed {\em globally}, i.e.,  for all $( \hat{p}(0),\col(\gamma_i(0)))\in\R^n\times [0, 1]^N $. The proposed mechanism is given by   
	\begin{equation}\label{soft_nudge}
		\dot{\hat{p}}(t)=\sum_{i\in\cali}x_i(t)-x^*+\frac{1}{\varepsilon}\left(\proj_{\calb}\left(\hat{p}(t)\right)-\hat{p}(t)\right),
	\end{equation} 
	where $ \calb $ is defined in \eqref{ball} and  $ \varepsilon>0 $ is a design parameter. We note that   the explicit expression of the projection of $ \hat{p}(t) $ on to the ball $ \calb $ is as follows\footnote{This can be verified by \cite[Prop. 2.1.3(b)]{bertsekas99nonlinear}.}:
	\begin{equation}\label{proj}
		\proj_{\calb}(\hat{p}(t))=\begin{cases}
			\hat{p}(t), & \text{if }\hat{p}(t)\in \calb,\\
			p_0+\frac{\bar{\delta}(\hat{p}(t)-p_0)}{\|\hat{p}(t)-p_0\|}, & \text{otherwise.}
		\end{cases}
	\end{equation}

	In the mechanism \eqref{soft_nudge}, the term $\sum_{i\in\cali}x_i(t)-x^*$ provides a suitable integral action as before to steer the aggregative behavior towards $x^*$. 
	However, different from \eqref{hard_nudge}, this term is outside the projection operator, and solutions of \eqref{soft_nudge} need not belong to the ball $\calb$ at all times. To emphasize this feature, we denote the dynamics \eqref{soft_nudge} as {\em soft}\, nudge\footnote{For related work on replacing projected dynamical systems with dynamics consisting of a penalty term, as in~\eqref{soft_nudge}, see the anti-windup approximation scheme studied in \cite{hauswirt2020robust}.}. 
	We note that outside the ball $\calb$, the term $\proj_{\calb}\left(\hat{p}(t)\right)-\hat{p}(t)\,$ is nonzero with the penalty gain $\varepsilon^{-1}$, thus attracting the price prediction $\hat p(t)$ to the ball and preventing the loss of trust. The parameter $\varepsilon$ is chosen sufficiently small such that trust variables  increase and reach the value of $1$ in finite time. 
	Below we establish the convergence properties of the  soft nudge mechanism.
	
	\begin{theorem}\label{thm:soft}
		Consider the closed-loop system formed by agents' model \eqref{agent} and the soft nudge mechanism  \eqref{soft_nudge} with   $ x^*\in \calx^* $. 
		Then, for any initial condition $\left(\hat{p}{(0)},\col (\gamma_i{(0)})\right)\in \R^n\times [0,1]^{N} $, there exists a  bounded Carath\'{e}odory solution  $t \mapsto \left(\hat{p}{(t)},\col (\gamma_i{(t)})\right)$ of the closed-loop system over the domain $ [0,\infty) $.
		Moreover, there exists some $ \varepsilon^*>0 $ such that for all $ \varepsilon\in(0, \varepsilon^*] $, any solution  $(\hat{p}{(t)},\col( \gamma_i{(t)}))$	
		converges to $ (p^*,\bone_N) $ with $ p^* $ given by \eqref{p_star}. Consequently, $\sum_{i\in \cali} x_i(t)$  converges to $x^*$ as desired. 
	\end{theorem} 
	\BP
	The proof is  divided  into  three parts. In the first part, we show that for any given $ (\hat{p}{(0)},\col( \gamma_i{(0)}))\in \R^n\times [0,1]^{N} $, there exists a bounded Carath\'{e}odory solution of \eqref{agent} and \eqref{soft_nudge}. The second part argues  that there exists some $ {\varepsilon}^*>0 $ such that  for all $ \varepsilon \in (0,{\varepsilon}^*] $,  the price prediction converges exponentially fast to the neighborhood of  the ball $ \calb $. We prove convergence of the solution to the  point $ (p^*,\bone_N) $ in the  last part.

	\textit{Existence  of  solutions:} By using \eqref{agent} and \eqref{soft_nudge}, we write the dynamics of the overall closed-loop system as
	\begin{subequations}\label{soft_cl}
		\begin{align}
			\dot{\hat p}&=h(\hat{p},\col(\gamma_i))\label{soft_cl_p} ,
			\\
			\dot\gamma_i&=\Pi_{[0,1]}\left(\gamma_i,\eta_i \psi_i(\|p(t)-\hat{p}\|)\right),\quad \forall i\in\cali\label{soft_cl_gam},
		\end{align}
	\end{subequations}
	where $ h(\hat{p},\col(\gamma_i)) := \sum_{i\in\cali}\pi_i({\hat{p}}, \gamma_i) - x^* + \frac{1}{\varepsilon}(\proj_{\calb}(\hat{p})-\hat{p}) $. Noting the nonexpansive property of $ \proj_{\calb} $ \cite[Prop. 2.1.3(c)]{bertsekas99nonlinear} and the definition of $ \pi_i $ given by \eqref{pi_def}, the map $ (\hat{p},\col(\gamma_i)) \mapsto h(\hat{p},\col(\gamma_i))$ is locally Lipschitz  in its arguments. Also, as discussed in the proof of Theorem~\ref{thm:hard}, we have that $ (\hat{p},t) \mapsto \psi_i(\|p(t)-\hat{p}\|) $ is  locally Lipschitz in $\hat{p}$ and measurable in $t$. Consequently,  existence of solutions follows  
	by showing that the hypotheses (i)-(iii) of Lemma~\ref{lem:exis2} are satisfied.

	We use the expression of $ \pi_i$ given by \eqref{pi_def} and rewrite the dynamics \eqref{soft_cl_p} as follows:
	\begin{equation}\label{p_mod_func_p}
		\begin{aligned}
			\dot{\hat{p}}=-\Big(\frac{1}{\varepsilon}I_n+\sum_{i\in\cali}\gamma_i  Q_i^{-1}\Big)\big(\hat{p}-\proj_{\calb}\left(\hat{p}\right)\big)+\nu ,
		\end{aligned}
	\end{equation}
	where $\nu:=\sum_{i\in\cali}(c_i+\gamma_{i}Q_i^{-1}(\hat\lambda_i-\proj_{\calb}(\hat{p})))-\sum_{i\in\cali}Q_i^{-1}\hat\lambda_i-x^* .$
	Note that the term $ \nu $ is bounded for all $ \hat{p}\in\R^n $ and $ \gamma_i\in [0,1] $. In particular, it follows from   $ \proj_{\calb}\left(\hat{p}\right)\in  \calb$  that there exists a  constant $ \bar{\nu}>0 $ such that  $ \|\nu\|\leq \bar{\nu}$ for all $ (\hat{p},\col(\gamma_{i}))\in \R^n \times [0,1]^N$.
	
	Now consider the  Lyapunov candidate
	$
	V\left(\hat{p}\right):=\frac{1}{2}\|\hat{p}-\proj_{\calb}\left(\hat{p}\right)\|^2 .
	$
	Since $ \proj_{\calb}\left(\hat{p}\right) $ is unique at any point $\hat p\in \R^n$ (cf. equation \eqref{proj}), it follows from   Danskin's Theorem \cite[Prop. B.25(a)]{bertsekas99nonlinear} that $ V(\hat{p}) $ is differentiable and $ \nabla V(\hat{p})=\hat{p}-\proj_{\calb}\left(\hat{p}\right) $.  Therefore  $ V(\hat{p}) $ satisfies Lemma \ref{lem:exis2}(i)-(ii). We next establish existence of solutions by  analyzing the inner product of  $ \nabla V(\hat{p}) $ and the right-hand side of \eqref{p_mod_func_p}. Recalling that $ h(\hat{p},\col(\gamma_i)) $ denotes the right-hand side of \eqref{p_mod_func_p} (cf. equation \eqref{soft_cl_p}), this inner product is computed as 
	\begin{multline*}
		\nabla V(\hat{p})^\top h(\hat{p},\col(\gamma_i))=-\|\hat{p}-\proj_{\calb}\left(\hat{p}\right)\|_{\sum_{i\in\cali}\gamma_{i}Q_i^{-1}}^2\\
		-\frac{1}{\varepsilon}\|\hat{p}-\proj_{\calb}\left(\hat{p}\right)\|^2+\left(\hat{p}-\proj_{\calb}\left(\hat{p}\right)\right)^\top \nu.
	\end{multline*}
	The first term on the right-hand side of the above equation is nonpositive as $ \gamma_{i}\in [0,1] $ and $ Q_i\succ 0 $ for all $ i\in \cali $. Using this fact and the bound on $\nu$, 
	we get
	\begin{multline}\label{in-prod-v-h}
		\nabla V(\hat{p})^\top h(\hat{p},\col(\gamma_i))\leq-\frac{1}{2\varepsilon}\|\hat{p}-\proj_{\calb}\left(\hat{p}\right)\|^2\\
		-\|\hat{p}-\proj_{\calb}\left(\hat{p}\right)\|\Big(\frac{1}{2\varepsilon}\|\hat{p}-\proj_{\calb}\left(\hat{p}\right)\| -\bar{\nu}\Big).
	\end{multline} 
	This implies that the inner product $ \nabla V(\cdot)^\top h(\cdot) $ is negative for all $ \|\hat{p}\|\geq \|p_0\|+ \bar\delta +2\varepsilon \bar{\nu} $ and $ \gamma_{i}\in [0,1] $. Therefore,  hypothesis (iii) of Lemma \ref{lem:exis2} is  satisfied, and the closed-loop system has a bounded Carath\'{e}odory solution for all $ t\geq 0 $.
	
	%
	%
	\noindent \textit{Convergence of $ \hat{p} $ to the neighborhood of $ \calb $:}	 Let a constant $ \tilde{\delta}>0 $ satisfying
	\begin{equation}\label{delta_tilde}
		\bar \delta<\tilde{\delta}< \min_{i\in\cali}\delta_i- \|\Delta p(t)\|\,,\quad \forall t\geq 0\,.
	\end{equation}
	Note that such $ \tilde{\delta} $ exists due to the condition \eqref{delta_bar}. Moreover, we deduce from \eqref{proj} that $ \|\hat{p}-\proj_{\calb}\left(\hat{p}\right)\|= \tilde{\delta}-\bar \delta$ for all $ \hat{p}\in \bd( \bar{B}(p_0,\tilde \delta))$. 
	Let $ \varepsilon \in (0,\varepsilon^*] $ with
	\begin{equation}\label{soft_ep}  
		\varepsilon^*:=\frac{\tilde{\delta}-\bar{\delta}}{2\bar{\nu}}.
	\end{equation}
	It then follows from \eqref{in-prod-v-h} that the time-derivative of the evolution of $V$ along any solution of~\eqref{soft_cl_p} satisfies 
	$
	\dot{V}\leq -\frac{1}{2\varepsilon}\|\hat{p}-\proj_{\calb}\left(\hat{p}\right)\|^2
	$
	for all  
	$\|\hat{p}-\proj_{\calb}\left(\hat{p}\right)\|\geq \tilde{\delta}-\bar{\delta}$.  
	Noting the definition of $ V $, we can then write  
	$
	\dot{V}\leq -\frac{1}{\varepsilon}V
	$
	whenever  
	$\|\hat{p}-\proj_{\calb}\left(\hat{p}\right)\|\geq \tilde{\delta}-\bar{\delta}$ or equivalently $ \hat{p}(t)  \not \in \bar{B}(p_0,\tilde \delta)$. 
	As a result, for
	any solution $t \mapsto (\hat{p}(t),\col(\gamma_i(t)))$, we have $ V(t)\leq V(0)\exp(-t/\varepsilon) $ as long as $\hat{p}(t)  \not \in \bar{B}(p_0,\tilde \delta)$.   
	Hence, if $\hat{p}$ is initialized outside the ball $ \bar{B}(p_0,\tilde \delta) $,  
	then it  converges exponentially fast to the ball in the time interval $ [0,T_1] $ with 
	\begin{equation}\label{t1}
		T_1=\varepsilon \ln \left(\frac{2V(0)}{(\tilde{\delta}-\bar{\delta})^2}\right),
	\end{equation}
	and we have $ \hat{p}(t)\in\bar{B}(p_0,\tilde \delta) $ for all $ t\geq T_1 $. Moreover, note that if $\hat p$ is initialized inside the ball $ \bar{B}(p_0,\tilde \delta) $, then it belongs to the ball for all $ t\geq T_1=0 $, since $\dot V$ is negative on 
	$\bd(\bar B(p_0, \tilde \delta))$.
	The above given reasoning establishes convergence of $ \hat{p}(t) $ to the ball $ \bar{B}(p_0,\tilde \delta) $ in finite time. 
	
	\emph{Convergence of $ (\hat{p},\col(\gamma_{i})) $ to $ (p^*,\bone_N) $:} For the rest of the proof we assume that $\varepsilon \in (0,\varepsilon^*]$ where $\varepsilon^*$ is given in~\eqref{soft_ep}. Consider any solution $t \mapsto (\hat{p}(t),\col(\gamma_i(t))$ of the closed-loop system. We divide the convergence analysis into three  time intervals $ [0,T_1] $, $ [T_1,T_2] $, and $ [T_2,\infty) $. Here, $T_1$ is equal to zero if $\hat p(0)\in \bar B(p_0, \tilde \delta)$, and $T_1$ is given by~\eqref{t1} otherwise. In other words, $T_1$ is the time when the trajectory $t \mapsto \hat{p}(t)$ enters and stays in the set $\bar{B}(p_0, \tilde \delta)$. Recall that $\gamma_i(t) \in [0,1]$ at all times. We will next show that full trust of all the agents is achieved in the time interval $[T_1,T_2]$ for some finite time $T_2$.

	Noting that $\tilde{\delta}$  satisfies \eqref{delta_tilde} and $\hat{p}(t)\in\bar B(p_0,\tilde \delta)$ in the time interval $ [T_1,\infty) $,  
	there exists some $\bar{\rho}>0$ such that $\|p(t)-\hat{p}(t)\|\leq \bar{\rho}<\min_{i\in\cali}\delta_i$ in the same time interval. By Assumption \ref{asm:psi}, we deduce that $ \psi_i(\|p(t)-\hat{p}(t)\|)\geq\psi_i(\bar\rho)>0 $ for all $ i\in \cali $.  This implies that, analogous to the discussions of trust variables in  the proof of Theorem~\ref{thm:hard} and \eqref{trust_time}, we have $ {\gamma}\ut{i}{t} = 1$ for all $t \geq T^i$, where $T^i :=  T_1+(1-\gamma\ut{i}{T_1})/(\eta_i\psi_i(\bar\rho)) $. Setting $ T_2:=\max_{i\in\cali}T^i $, we conclude that   $ \col( \gamma_i{(t)})=\bone_N $ for all $ t\in [T_2,\infty) $, i.e., full trust of the agents is obtained in the time interval $[T_1,T_2]$. 
	
	In the time interval $ [T_2,\infty) $, using $ {\gamma}\ut{i}{t}=1 $ for all $ i\in\cali $, the dynamics of the price prediction~\eqref{p_mod_func_p} reduces to  
	\begin{equation}\label{eq:phat-t2}
		\dot{\hat{p}}=-\sum_{i\in\cali}Q_i^{-1}\left(\hat{p}-{p}^*\right)+\frac{1}{\varepsilon}\left(\proj_{\calb}\left(\hat{p}\right)-\hat{p}\right),
	\end{equation}
	where $ \hat{p}{(T_2)}\in \bar{B}(p_0,\tilde \delta) $ and we used the expression of $ x^* $ in \eqref{x_star}. Now, we consider the  Lyapunov  candidate $ W\left(\hat{p}\right):=\frac{1}{2}\|\hat{p}-p^*\|^2 $ and analyze its evolution along the solution of~\eqref{eq:phat-t2}. We have  
	\begin{align*}
		\dot{W}=-\|\hat{p}-{p}^*\|_{\sum_{i\in\cali}Q_i^{-1}}^2+\frac{1}{\varepsilon}\left(\hat{p}-{p}^*\right)^\top(\proj_{\calb}(\hat{p})-\hat{p}).
	\end{align*}
	The second term on the right-hand side  satisfies
	\begin{equation}\label{p-tilde-ineq}
		\begin{aligned}
			&\left(\hat{p}-{p}^*\right)^\top\left(\proj_{\calb}\left(\hat{p}\right)-\hat{p}\right)\\
			&=\left(\hat{p}-\proj_{\calb}\left(\hat{p}\right)\right)^\top\left(\proj_{\calb}\left(\hat{p}\right)-\hat{p}\right)\\
			&+\left(\proj_{\calb}\left(\hat{p}\right)-{p}^*\right)^\top\left(\proj_{\calb}\left(\hat{p}\right)-\hat{p}\right)\leq 0,
		\end{aligned}
	\end{equation}
	where we used $ p^*\in \calb $ and \cite[Prop. 2.1.3(b)]{bertsekas99nonlinear} to write the inequality. Consequently, we obtain 
	$
	\dot{W}\leq-\left(\hat{p}-{p}^*\right)^\top\sum_{i\in\cali}Q_i^{-1}\left(\hat{p}-{p}^*\right).
	$
	This implies that $ \hat{p} $ exponentially  converges to $ p^* $ in the time interval $[T_2,\infty)$, 
	and the aggregate behavior $\sum_{i\in \cali} x_i$ converges to $x^*$. 
	\EP
	
	\begin{remark}
		While Theorem \ref{thm:soft} guarantees existence of a sufficiently small $\varepsilon^*$ given by \eqref{soft_ep}, computing its value requires the knowledge of bounds on agent parameters $ c_i $, $Q_i$, $\delta_i$, and $\hat\lambda_i$. If such bounds are not available, one can opt for the hard nudge mechanism \eqref{hard_nudge} at the cost of restricting the initial condition $\hat p(0)$ to $\mathcal{B}$. 
		\oprocend
	\end{remark}
		\begin{remark}\label{rem:convex-cost}
						The results of the hard and soft nudge mechanisms remain valid for more general classes of cost functions than \eqref{cost}. In particular,  let  the cost functions be of the form $
			J_i(z_i,p):=c_i(z_i)+z_i^\top p
			$, 	where  $ c_i:\R^n\to \R $ is $ \calc^2 $ and strongly convex. It follows that the model of the agents in \eqref{agent} will be modified to
			\begin{equation*}
			\agent_i:\left\{\begin{split}
			\dot\gamma_i(t)&=\Pi_{[0,1]}(\gamma_i(t),\eta_i \psi_i(\|p(t)-\hat{p}{(t)}\|)),\\
			x_{i}(t)&=(\nabla c_i)^{-1}\left(-\gamma_i(t) \, \hat{p}(t)-\left(1-\gamma_i(t)\right)\hat\lambda_i\right).
			\end{split}\right.
			\end{equation*}
			It can be shown that for any desired behavior $ x^*\in \calx^* $ with
			\begin{equation*}
			\calx^*:=\Big\{x\in\R^n\mid x=\sum_{i\in\cali}(\nabla c_i)^{-1}(-\hat{p}),\,   \hat{p}\in \calb\Big\},
			\end{equation*}
			both hard and soft nudges guarantee convergence of the aggregative behavior to $ x^*$. However, when the desired behavior is time-varying, as considered in the next section,  devising a suitable nudge mechanism becomes much more challenging. Therefore, to unify the presentation throughout the paper, we have provided our results for the linear-quadratic cost function \eqref{cost}. \oprocend
	\end{remark}

	\section{A nudge mechanism for temporal desired behaviors}\label{s:adaptive}
	So far, we have treated the desired aggregative behavior as a fixed point. However, this point may vary with time in practice due to changes in the market condition, the climate, and government policies.  
	In the context of power systems, for instance, climate change affects the efficiency of power production as well as the energy consumption \cite{climate2014us}. The policies passed  by the government also affect the market substantially, see e.g. \cite{yin2010state} regarding renewable energy. These changes entail variations of the desired aggregative behavior over time. Building on \eqref{soft_nudge}, we design here a nudge mechanism that steers the aggregative behavior of the agents to a desired time-varying 
	signal $t\mapsto x^*(t) $. 
	The set of admissible reference signals $x^*(\,\cdot\,)$ is given by the assumption below. 
	\begin{assumption}\label{asm:xs-t}
		The signal $t\mapsto x^*(t) $  belongs to the set $ \calx^* $ given by \eqref{xs_def} for all $ t\in[0,\infty) $. In addition,  $x^*(\,\cdot\,)$ is continuously 		differentiable with bounded
		derivative over the domain $[0,\infty)$, that is, there exists a constant $\theta > 0$ such that $ \|\dot{x}^*(t)\|\leq \theta $ for all $ t\in[0,\infty) $. \oprocend 
	\end{assumption}
	The above assumption indicates that the desired aggregative behavior of the agents satisfies a regularity condition in the sense that it is smooth and belongs to the admissible set $ \mathcal{X}^* $. For all $t \in[0,\infty)$, since $x^*(t) \in \mathcal{X}^*$, we obtain from \eqref{xs_def} that there exists a unique $ p^*(t)\in \calb$  
	such that 
	\begin{equation}\label{xs-ps-t}
		x^*(t)=\sum_{i\in\cali}\left(c_i-Q_i^{-1}{p}^*(t)\right).  
	\end{equation}
	Rearranging the terms, $ p^*(t) $ can be written explicitly as  
	\begin{equation}\label{ps-t}
		p^*(t)=\Big(\sum_{i\in\cali}  Q_i^{-1} \Big) ^{-1} \big(-x^*(t)+ \sum_{i\in\cali} c_i\big). 
	\end{equation}
	Note from Assumption~\ref{asm:xs-t} that the signal $t \mapsto p^*(t)$ 
	is differentiable with a bounded derivative. 
	If the system regulator had accurate knowledge of all $ Q_i$ and $ c_i $ parameters, it could have obtained the desired behavior by setting the price prediction equal to $ p^*(t) $. However, since the cost functions of the agents are unknown to the system designer, such a simple strategy cannot be implemented. This asks for a more sophisticated design, and to that end, we propose the following \textit{adaptive} nudge mechanism 
	\begin{subequations}\label{adap_nudge}
		\begin{align}
			\phantom{\dot{K}(t)}
			&\begin{alignedat}{2}\label{adap_nudge_one}
				\mathllap{\dot{\hat{p}}(t)} &= \sum_{i\in\cali}x_i(t)&&-x^*(t)+K(t)\dot{x}^*(t)\\[-5pt]
				& &&+\frac{1}{\varepsilon}\left(\proj_{\calb}\big(\hat{p}(t)\right)-\hat{p}(t)\big),
			\end{alignedat} \\
			&\begin{alignedat}{2}\label{adap_nudge_two}
				\mathllap{\dot{K}(t)} &= \tau\Big(\sum_{i\in\cali}x_i(t)&&-x^*(t)\Big){\dot{x}^*(t)}^\top\\[-5pt]
				& &&-\tau\,\sigma_s\big(\|K(t)\|_\text{F}\big) K(t),
			\end{alignedat} 
		\end{align}
	\end{subequations}
	where $ \calb $ is given by \eqref{ball}, $ \|K(t)\|_\text{F} $ is the Frobenius norm of $ K(t) $, $ \varepsilon>0 $, $\tau>0$, and the function $\sigma_s:\R_{\ge 0}\to [0,\sigma]  $ is given by
	\begin{align}\label{sigma}
		\sigma_s(u) :=\begin{cases}
			0 & \text{ if }u<k_0, \\
			\sigma\left(\frac{u}{k_0}-1\right) & \text{ if }k_0\leq u\leq 2k_0, \\
			\sigma & \text{ if } 2k_0<u\,.
		\end{cases}
	\end{align}
	In the above definition, $ \sigma>0 $ and $ k_0>0 $ are design parameters that are selected afterwards.  	
	
	\noindent 	{\em Interpretation of the adaptive nudge mechanism:} There are several remarks in order concerning the adaptive nudge \eqref{adap_nudge}: (i) This mechanism simplifies to the soft nudge mechanism \eqref{soft_nudge} in case of a stationary desired aggregative behavior. Namely, with $\dot x^*(t) = 0$, the dynamics \eqref{adap_nudge_one} reduces to \eqref{soft_nudge} and \eqref{adap_nudge_two} can be discarded. (ii) Compared to the soft nudge mechanism, the additional term $K(t) \dot x^*(t)$ is included to cope with the temporal nature of the desired aggregative behavior by tracking the signal $\dot p^*(t)$ given by (cf. equation \eqref{ps-t}) 
	\begin{equation}\label{e:p-star-dot}
		\dot p^*(t)= K^* \dot x^*(t), \quad K^*:=-\Big(\sum_{i\in\cali}Q_i^{-1}\Big)^{-1}.
	\end{equation}
	Again since the regulator is not aware of all cost functions, a static choice $K(t)=K^*$  
	would not be feasible and we, therefore, appeal to the adaptive law \eqref{adap_nudge_two}.  
	(iii) The first term on the right-hand side of  \eqref{adap_nudge_two} is chosen such that sign-indefinite terms in the time-derivative of the Lyapunov function are canceled out.    
	The second term provides a state-dependent damping that prevents 
	the matrix $K(t)$ to become unbounded.

	\noindent \textit{Selection of design parameters:} In order to guarantee convergence of the adaptive nudge algorithm, the design parameters $ \varepsilon $, $ \sigma$, and $ k_0$ should be chosen appropriately.  The treatment in Lemma~\ref{lem:adap-uub} in the appendix suggests to choose $ \varepsilon\in\cali_{\varepsilon} $, $ \sigma\in\cali_{\sigma} $, and $ k_0\in\cali_{k_0} $  with
	\begin{equation}\label{e:adap-par}
		\begin{split}
			\cali_{\varepsilon}&:=\Big(0,\theta^{-1}(1+\lambda_{\max}(\sum_{i\in\cali}Q_i^{-1}))^{-1}\Big],\\
			\cali_{\sigma}&:=\Big[ 2\theta(1+\lambda_{\max}(\sum_{i\in\cali}Q_i^{-1})),\infty\Big),\\
			\cali_{k_0}&:=\Big[\sqrt{n}\lambda_{\max}(\sum_{i\in\cali}Q_i^{-1})/\lambda_{\min}^2(\sum_{i\in\cali}Q_i^{-1}),\infty\Big).
		\end{split}
	\end{equation} 
	Note that the design parameters can take any values within the bounds indicated above, and therefore their selection is oblivious of the exact values of the cost parameters. 
	
	The main result of this section is provided in the following theorem.

	\begin{theorem}\label{t:new-adaptive}
		Consider the closed-loop system formed by agents' model \eqref{agent} and the adaptive nudge mechanism  \eqref{adap_nudge} with   $t\mapsto x^*(t)$ satisfying Assumption~\ref{asm:xs-t}. Let the design parameters satisfy $ \sigma\in\cali_{\sigma} $ and $ k_0\in\cali_{k_0} $ with the intervals $ \cali_{\sigma} $ and $ \cali_{k_0} $ given by \eqref{e:adap-par}. Then, there exists some $ \varepsilon^*\in\cali_{\varepsilon} $ with $ \cali_{\varepsilon} $ given by \eqref{e:adap-par}  such that for all $ \varepsilon\in (0,\varepsilon^*] $ and any initial condition $ (\hat{p}{(0)},K(0),\col( \gamma_i{(0)}))\in \R^n\times \R^{n\times n} \times [0,1]^{N}$, there exists a bounded Carath\'{e}odory solution  $t \mapsto \left(\hat{p}(t),K(t),\col (\gamma_i(t))\right)$ of the closed-loop system over the domain $ [0,\infty) $. Moreover, any solution   
		$ (\hat{p}{(t)},\col( \gamma_i{(t)})) $ converges to $ (p^*(t),\bone_N) $ with $ p^*(t) $ given by \eqref{ps-t}. Consequently,  $ \sum_{i\in \cali}x_i(t) $ converges to $x^*(t)$ as desired.
	\end{theorem}
	\BP
	Our proof builds on the results of Lemma~\ref{lem:adap-uub}. Let $ \varepsilon\in\cali_{\varepsilon} $, $ \sigma\in\cali_{\sigma} $, and $ k_0\in\cali_{k_0} $, then it follows from Lemma~\ref{lem:adap-uub} that the closed-loop system admits a bounded  Carath\'{e}odory solution over domain $ [0,\infty) $. Consider any solution $t \mapsto \left(\hat{p}(t),K(t),\col (\gamma_i(t))\right)$.  
	Again from Lemma~\ref{lem:adap-uub}, there is a finite time $ T\geq 0 $ such that for all $ t\geq T $, we have $ \|\hat{p}(t)\|\leq \bar p $ and $ \|K(t)\|_{\text{ F }}\leq \bar k $ with $ \bar p $ and $ \bar k $ given by  \eqref{uub}.  Next we prove convergence of $ (\hat{p}{(t)},\col( \gamma_i{(t)})) $ to $ (p^*(t),\bone_N) $ by considering three time intervals $ [T,T_1] $, $ [T_1,T_2] $, and $ [T_2,\infty) $. The first time interval concerns the convergence analysis of $ \hat{p}(t) $  to the neighborhood of $ \calb $. Full trust of the agents is achieved in the second time interval, while convergence of $\hat{p}(t) $ to $ p^*(t)$ is established in the last time interval.

	We analyze the  interval $ [T,T_1] $ by considering the price prediction dynamics \eqref{adap_nudge_one} as a system with bounded exogenous signals. In particular, we substitute  the expression of $ x_i $ given by \eqref{opt_action_cont} and \eqref{pi_def} into 
	\eqref{adap_nudge_one} to get:
	\begin{equation*}
		\begin{aligned}
			\dot{\hat{p}}&=-\Big(\frac{1}{\varepsilon}I_n+\sum_{i\in\cali}\gamma_i(t)  Q_i^{-1}\Big)\big(\hat{p}-\proj_{\calb}(\hat{p})\big)+\nu(t),
		\end{aligned}
	\end{equation*}
	where $ t\mapsto\gamma_{i}(t) $ and $t\mapsto \nu(t) $ are treated as exogenous signals and $\nu(t):=\sum_{i\in\cali}(c_i+\gamma_{i}(t)Q_i^{-1}(\hat\lambda_i-\proj_{\calb}(\hat{p})))-\sum_{i\in\cali}Q_i^{-1}\hat\lambda_i-x^*(t)+K(t)\dot{x}^*(t) .$
	From  the proof of Lemma~\ref{lem:adap-uub}, we see that the time instant $ T $ and the ultimate  bounds $ \bar p $ and $ \bar k $ are uniform for all $ \varepsilon\in\cali_{\varepsilon} $.
	This, in addition to $\proj_{\calb}\left(\hat{p}\right)\in  \calb$, $ \gamma_i(t)\in [0,1] $, and boundedness of $ x^*(t) $ and $  \dot{x}^*(t) $ (cf. Assumption \ref{asm:xs-t}), imply that  $\nu(t) $ is  uniformly ultimately bounded. More precisely, there exists some constant $ \bar \nu>0 $ such that $ \|\nu(t)\|\leq \bar \nu $ for all $ t\geq T $ and all $ \varepsilon\in\cali_{\varepsilon} $. Next we use this property and show that suitable selection of $ \varepsilon $ provides convergence of $ \hat{p}(t) $  to the neighborhood of $ \calb $  in finite time.
	Let 
	\begin{equation*}
		\varepsilon^*:=\min\Big\{\frac{\tilde{\delta}-\bar{\delta}}{2\bar{\nu}},\,\theta^{-1}(1+\lambda_{\max}(\sum_{i\in\cali}Q_i^{-1}))^{-1}\Big\},
	\end{equation*}
	with $ \tilde{\delta} $ satisfying \eqref{delta_tilde}. This results in $ \varepsilon^*\in\cali_{\varepsilon} $. Moreover, following the steps of the  proof of Theorem~\ref{thm:soft}, there exists some  $ T_1\geq T $ such that by choosing $ 0<\varepsilon\leq \varepsilon^* $, $ \hat{p}(t) $  belongs the ball $ \bar{B}(p_0,\tilde \delta) $ for all $ t\geq T_1 $. We note that such selection of $ \varepsilon $ is possible since $ \bar\nu $, and hence $ \varepsilon^* $, are independent of the choice of $ \varepsilon\in\cali_{\varepsilon}$. 
	
	Bearing in mind $ \hat{p}(t)\in \bar{B}(p_0,\tilde \delta)$ for all $ t\geq T_1 $, an analogous argument to  the proof of  Theorem~\ref{thm:soft} can be used to show that  
	there exists  a finite time $ T_2\geq T_1 $ such that we have $ \gamma_{i}(t)=1 $ for all $ i\in\cali $ and $ t\geq T_2 $. Next we exploit $ \gamma_{i}(t)=1 $ to establish convergence of $ \hat{p} $ to $ p^* $ in the time interval $[T_2,\infty) $. We perform a change of coordinates to ease the notation, namely, $ (\hat{p},K)\mapsto (\tilde{p},\Phi) $ with $ \tilde{p}=\hat{p}-p^* $ and $  \Phi=K-K^* $ 
	where $ K^* $ is given by \eqref{e:p-star-dot}.  In these coordinates, the closed-loop system, comprised of  \eqref{agent} and \eqref{adap_nudge},  takes the form
	\begin{equation}\label{track-dyn-last}
		\begin{aligned}
			\dot{\tilde{p}}&=-\sum_{i\in\cali}Q_i^{-1}\tilde{p}+\Phi\,\dot{x}^*(t)+\frac{1}{\varepsilon}\left(\proj_{\calb}\left(\hat{p}\right)-\hat p\right),\\
			\dot{\Phi}&=-\tau\sum_{i\in\cali}Q_i^{-1}\tilde{p}\,\dot{x}^{*}(t)^\top-\tau \underbrace{\sigma_s\big(\|\Phi+K^*\|_\text{F}\big) (\Phi+K^*)\,}_{\sigma_s(\|K\|_\text{F}) K},	
		\end{aligned}
	\end{equation}
	where we have used $ \gamma_{i}(t)=1 $ and the  expressions of $ \pi_i $, $ x^*(t) $, and $ \dot{p}^*(t) $, respectively given by \eqref{pi_def}, \eqref{xs-ps-t}, and \eqref{e:p-star-dot}. For the rest of the proof, we use the following definition for notational simplicity. 
	\begin{equation}\label{q-def}
		Q:=\sum_{i\in\cali}Q_i^{-1}\,.
	\end{equation}
	Consider the following Lyapunov  candidate
	$
	V(\tilde{p},\Phi):=\frac{1}{2}\|\tilde{p}\|^2+\frac{1}{2\tau}\Tr\left(\Phi^\top Q^{-1}  \Phi\right).
	$
	The evolution of $V$ along the solutions of ~\eqref{track-dyn-last} is given by  
	\begin{align*}
		\dot{V}=&-\|\tilde{p}\|_Q^2+\tilde{p}^\top\Phi\,\dot{x}^*(t)+\frac{1}{\varepsilon}\tilde{p}^\top\left(\proj_{\calb}\left(\hat{p}\right)-\hat p\right)\\
		&-\Tr\left(\dot{x}^{*}(t) \,\tilde{p}^\top \Phi \right)-\sigma_s\big(\|K\|_\text{F}\big)\Tr\left(K^\top Q^{-1}  \Phi\right).
	\end{align*}
	It follows from $ \tilde{p}^\top\Phi\,\dot{x}^*(t)=\Tr\left(\dot{x}^{*}(t)\, \tilde{p}^\top \Phi \right) $ and \eqref{p-tilde-ineq} that
	\begin{equation}\label{v2-dot}
		\dot{V}\leq-\|\tilde{p}\|_Q^2-\sigma_s\big(\|K\|_\text{F}\big)\Tr\left(K^\top Q^{-1}  \Phi\right).
	\end{equation}
	We proceed to show that, given $ k_0\in \cali_{k_0} $, the second term on the right-hand side is nonpositive. 
	We note that $ \Phi=K+Q^{-1} $ due to \eqref{e:p-star-dot} and \eqref{q-def}.  
	It then follows from $ \sigma_s(\,\cdot\,)\geq 0 $ that		
	\begin{multline}\label{sig_tra_up}
		-\sigma_s\big(\|K\|_\text{F}\big)\Tr\left(K^\top Q^{-1}  \Phi\right) 
		\leq-\frac{\sigma_s\big(\|K\|_\text{F}\big)}{\lambda_{\max}\left(Q\right)}\|K\|_\text{F}^2\\
		+\sigma_s\big(\|K\|_\text{F}\big)\|K\|_\text{F}\|Q^{-2}\|_\text{F}.
	\end{multline}
	In the previous inequality, we used  $ \Tr\left(K^\top Q^{-1}  K\right)\geq  \lambda_{\min}\left(Q^{-1}\right)\|K\|_\text{F}^2  $ and  $ \lambda_{\min}\left(Q^{-1}\right)=1/\lambda_{\max}\left(Q\right) $ to find the first term on the right-hand side, and the second term is obtained using  Cauchy–Schwarz inequality as $ |\!\Tr(K^\top Q^{-2})| \leq \|K\|_\text{F} \|Q^{-2}\|_\text{F}$. In addition, notice that we have $ \|Q^{-2}\|_\text{F}\leq \sqrt{n}/\lambda_{\min}^2(Q)  $. It then follows from the definition of  $ \cali_{k_0} $ that   
	$
	||Q^{-2}\|_\text{F}\leq {k_0}/{ \lambda_{\max}(Q)} 
	$ for all $ k_0\in \cali_{k_0} $.
	The latter implication implies that \eqref{sig_tra_up} can be further bounded as
	\begin{multline*}
		-\sigma_s\big(\|K\|_\text{F}\big)\Tr\left(K^\top Q^{-1}  \Phi\right) 
		\\ \leq -\frac{\sigma_s\big(\|K\|_\text{F}\big)}{\lambda_{\max}\left(Q\right)}\|K\|_\text{F}\big(\|K\|_\text{F}-k_0\big).
	\end{multline*}
	Bearing in mind the definition of $ \sigma_s(\,\cdot\,)\geq 0 $ given by  \eqref{sigma}, we find that $ \sigma_s\big(\|K\|_\text{F}\big)\left(\|K\|_\text{F}-k_0\right)\geq 0 $ for all $ K\in\R^{n\times n} $. Combining this with the above inequality results in $ -\sigma_s\big(\|K\|_\text{F}\big)\Tr\left(K^\top Q^{-1}  \Phi\right)\leq 0 $. Consequently, the relation  \eqref{v2-dot} provides
	\begin{equation}\label{adap-v-last}
		\dot{V}\leq-\|\tilde{p}\|_{Q}^2\,.
	\end{equation} 
	Next, recalling that the dynamics \eqref{track-dyn-last} is a nonautonomous system, we use Barbalat's lemma \cite[Lem. 4.2]{slotine1991applied} to conclude convergence of $ \tilde{p}(t) $ to the origin. Let $ f(t):=\int_{T_2}^{t} \|\tilde{p}(s)\|_{Q}^2 ds$ for $ t\geq T_2 $. From \eqref{track-dyn-last}, we see that $ \dot{\tilde{p}}(t) $ is bounded for all $ t\geq T_2 $. This implies that $ \ddot{f}(t) $ is bounded too, and thus $ \dot f(t) $ is uniformly continuous. The next step is to show that the function $ f(t) $ has a finite limit as $ t\to \infty $. For that, we integrate both sides of  \eqref{adap-v-last} and use the definition of $ f(t) $ with $ V(t)\geq 0 $ to obtain $$ \lim_{t\to \infty}f(t)\leq V(T_2). $$ The left-hand side of the inequality above is bounded since $ V(T_2) $ is bounded. 
	It then follows from Barbalat's lemma that $ \lim_{t\to\infty}\dot{f}(t)=0 $, i.e., 
	$  \tilde{p}(t)\to 0  $ as $ t\to \infty $. We conclude that $ \hat{p}(t) $ converges to $ p^*(t) $ in the time interval $ [T_2,\infty) $, and in turn, the aggregative behavior  $ \sum_{i\in \cali}x_i(t) $ converges to $x^*(t)$ as desired.
	\EP
	\begin{remark}\label{rem:adap-hard}
			We note that one can also devise an adaptive nudge mechanism that is built on the hard nudge \eqref{hard_nudge} as follows:
			\begin{subequations}
				\begin{align*}
					\phantom{\dot{K}(t)}
					&\begin{alignedat}{1}
						\mathllap{\dot{\hat{p}}(t)} &= \Pi_\calb\Big(\hat{p}(t),\sum_{i\in\cali}x_i(t)-x^*(t)+K(t)\dot{x}^*(t)\Big),
					\end{alignedat} \\
					&\begin{alignedat}{2}
						\mathllap{\dot{K}(t)} &= \tau\Big(\sum_{i\in\cali}x_i(t)&&-x^*(t)\Big){\dot{x}^*(t)}^\top\\[-5pt]
						& &&-\tau\,\sigma_s\big(\|K(t)\|_\text{F}\big) K(t),
					\end{alignedat} 
				\end{align*}
			\end{subequations}
			where  $ \|K(t)\|_\text{F} $ is the Frobenius norm of $ K(t) $, $ \tau>0 $, and the function $ \sigma_s:\R_{\geq 0}\to [0,\sigma] $ is defined in \eqref{sigma}. We can then show that for any $t\mapsto x^*(t)$ satisfying Assumption~\ref{asm:xs-t}, choosing the design parameters $ \sigma>0 $ and $ k_0\in\cali_{k_0} $ with $ \cali_{k_0} $ given by \eqref{e:adap-par}, results in convergence of the aggregative behavior  to  $  x^*(t) $.  We note, however, that the resulting convergence is restricted to the ball $\calb$ and is thus not global, unlike in the  adaptive (soft) nudge mechanism \eqref{adap_nudge}. The details of the analysis are omitted due to lack of space. \oprocend
	\end{remark}
	\section{Case study}\label{sec:sim}
	We illustrate the performance of our nudge mechanisms by considering the problem of coordinated charging  of plug-in electric vehicles \cite{ma2015distributed}.   In this problem, the objective of the regulator is to control the aggregative power demand over a charging horizon. 
	
	We consider a population of $ \cali=\{1,\dots,10\}  $ agents, where each agent  $ i $ aims at choosing its charging strategy over the charging horizon of length $ n=24 $, namely  $ z_i\in \calx_i\subset \R^n  $, such that its cost function given below is minimized:
	\begin{equation}\label{pev-cost}
		C_i(z_i,p):=a_i z_i^\top z_i + b_i z_i^\top \bone_n+z_i^\top p,
	\end{equation}
	where $ a_i\in[0.004,0.006] $ and $ b_i\in [0.065,0.085] $. The set $ \calx_i $ is nonempty, compact, and convex, and it is defined as follows:
	\begin{equation*}
		\calx_i:=\left\{z_i \in \R^n \mid  z_i\in [0,\bar x_i]^n,\, \bone_n^\top z_i= d_i\right\},
	\end{equation*}
	where 
	$ \bar x_i\in [8,10]$(kW) is the maximum charging rate at any instant, and $ d_i\in[25,35] $(kWh) is the total energy required by the agent.   
	
	Since agents choose their actions from the sets $ \calx_i $, rather than $\R^n$, the expression of the optimal action  \eqref{opt_action} modifies to \cite[Prop. 2.1.2 and 2.1.3(b)]{bertsekas99nonlinear}, 
	\begin{equation}\label{sim:const_act}
		x_i=\proj_{\calx_i}\Big(-\frac{1}{2a_i}\big(b_i \bone_n+\gamma_i\hat{p}+(1-\gamma_i)\hat{\lambda}_i\big)\Big).
	\end{equation}
	Note that for $\mathcal{X}_i=\R^n$, the expression \eqref{sim:const_act} reduces to \eqref{opt_action}. 
	As for the choice of  $\psi_i$, we pick  
	$ \psi_i(\|p-\hat{p}\|)=-\tanh(h_i(\|p-\hat{p}\|-\delta_i))  $ with $ h_i\in [2,5] $,  which satisfies Assumption \ref{asm:psi}, and we select $ \delta_i\in [0.3,0.5]$(\$/kWh), $ \eta_i\in[3,5] $,   $ \hat{\lambda}_i\in[0.1,0.5]^n $(\$/kWh),   $ \gamma_i(0)\in[0,0.7] $ to simulate the model.

	Taking Assumption \ref{std_asm} regarding the actual price signal into consideration, we pick $ p_0=0.3 \bone_n$(\$/kWh) and consider price fluctuations to satisfy $ \|\Delta p(t)\|\leq 0.1 $(\$/kWh) for all $ t\geq 0 $. 
	Let $ \rho=0.2 $, then $ \rho $ is less than or equal to the expression on the right hand side of \eqref{delta_bar}. Consequently, the open ball ${B}{(p_0,\rho)}=\{\hat{p}\in\R^n \mid \|\hat{p}-p_0\|< \rho\}$ is  a feasible set 
	for the price prediction such that the regulator can gain
	agents' trust. We also define the ball $ \calb $ by choosing $ \bar{\delta}=0.15 $. Therefore the condition \eqref{delta_bar} is satisfied noting that $\bar{\delta}< \rho  $.
	\subsection{Stationary desired behavior}
\begin{figure}
		\centering
		\includegraphics[width=0.8\columnwidth]{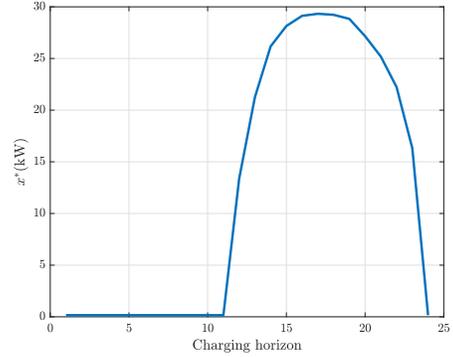}
		\caption{Desired stationary aggregative power demand over the charging horizon.}
		\label{fig:xstar-sta}
	\end{figure}

	Here we demonstrate convergence of the aggregative behavior to a desired behavior $x^*$ shown in  Fig. \ref{fig:xstar-sta},  
	under both hard and soft nudge mechanisms.   
	The desired aggregative behavior 
	specifies the goal of the system regulator in nudging the vehicles to charge their batteries in a specific interval. 
	
\begin{figure}
		\centering
		\includegraphics[width=\columnwidth]{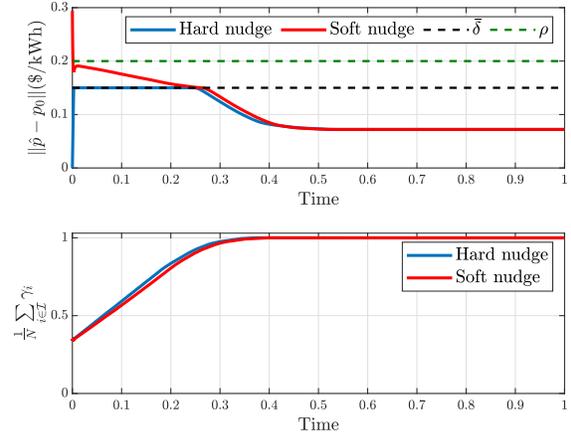}
		\caption{Distance of hard and soft nudges' price predictions to $ p_0 $ and the average of the trust variables.}
		\label{fig:p-gamma-sta}
	\end{figure}

\begin{figure}
		\centering
		\includegraphics[width=\columnwidth]{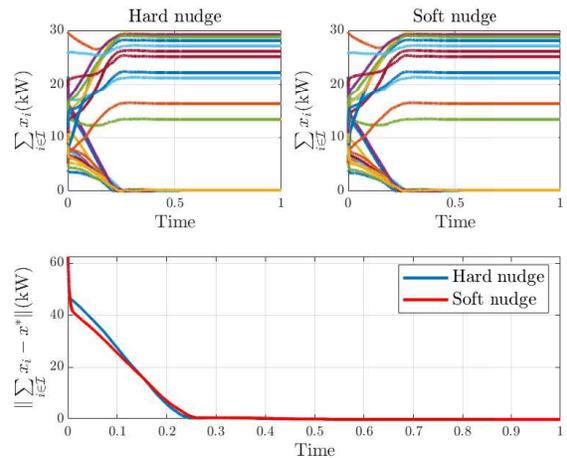}
		\caption{Aggregative power demands due to hard and soft nudges and their distance to the desired stationary power demand.}
		\label{fig:x-sum-sta}
	\end{figure}

	We choose $ \hat{p}(0)=p_0 \in \mathcal{B}$ for the hard nudge, whereas we  
	set  $ \varepsilon=10^{-3} $ and $ \hat{p}(0)=p_0+0.06\bone_n\notin \calb $ for the soft nudge to demonstrate convergence for an initialization outside the ball $\mathcal{B}$.  
	Fig. \ref{fig:p-gamma-sta} shows the distance of the mechanisms' price predictions to $ p_0 $  and the average of the trust variables. We observe that for the hard nudge, the price prediction belongs to the ball $ \calb $ for all times, and as a result, the trust variables converge to one. The latter is deduced from convergence of the average of the trust variables to one and $ \gamma_{i}\in[0,1] $. For the soft nudge, the price prediction converges to a positively invariant set inside the open ball $ {B}{(p_0,\rho)} $, which in turn increases the agents' trust on $ \hat{p} $. After gaining full trust of the agents,  the price predictions of both mechanisms converge to $ p^*\in\calb $. Therefore, the aggregative behavior of the agents, namely the aggregative power demand, converges to $ x^* $ as demonstrated in Fig. \ref{fig:x-sum-sta}.

	\subsection{Temporal desired behavior}	
\begin{figure}
		\centering
		\includegraphics[width=\columnwidth]{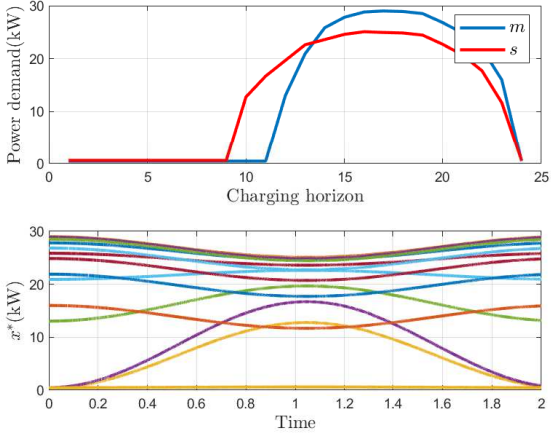}
		\caption{Evolution of the desired temporal aggregative power demand $ x^*(t)=\frac{1+\cos (3t)}{2}m+\frac{1-\cos (3t)}{2}s $.}
		\label{fig:xstar-tem}
	\end{figure}

\begin{figure}
		\centering
		\includegraphics[width=\columnwidth]{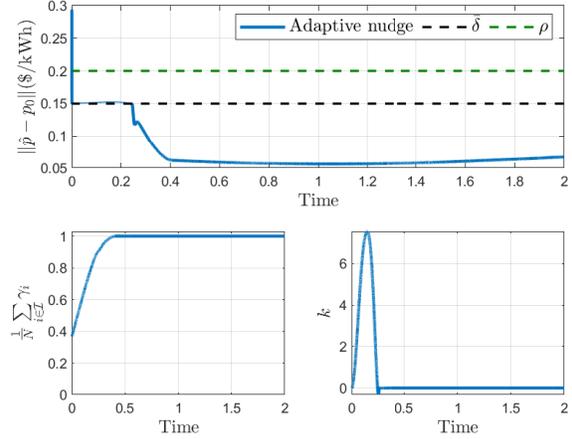}
		\caption{Distance of adaptive nudge's price prediction to $ p_0 $, the average of the trust variables, and evolution of the adaptive gain $ k $.}
		\label{fig:p-gamma-k-tem}
	\end{figure}

\begin{figure}
		\centering
		\includegraphics[width=\columnwidth]{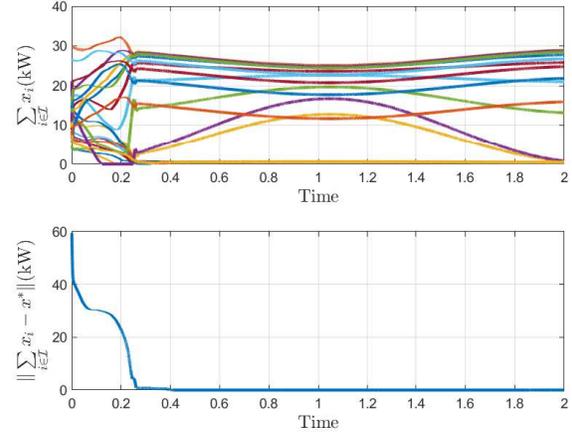}
		\caption{Aggregative power demands due to adaptive nudge and its distance to the desired temporal power demand.}
		\label{fig:x-sum-tem}
	\end{figure}

	Next, we consider the case  where the desired aggregative behavior varies with time, and employ the adaptive nudge protocol to steer the aggregative behavior towards such behavior. We choose  
	the desired behavior as $ x^*(t)=\frac{1+\cos (3t)}{2}m+\frac{1-\cos (3t)}{2}s $ with $ m $ and $ s $ shown in  Fig. \ref{fig:xstar-tem}.  
	Recalling the structure of the cost function $ C_i $ as \eqref{pev-cost}, we observe that its minimization is equivalent to minimization of $ J_i $ given by \eqref{cost} with $ Q_i=2 a_i I_n $ and $ c_i=-\frac{b_i}{2 a_i}\bone_n $. 
	Therefore, the matrix $ K^* $ in \eqref{e:p-star-dot} and thus the matrix $K$ in \eqref{adap_nudge} becomes a scalar matrix, i.e., $K=kI_n$, and the adaptive nudge \eqref{adap_nudge} reduces to  
	\begin{align*}
		\dot{\hat{p}}&=\sum_{i\in\cali}x_i-x^*(t)+k\,\dot{x}^*(t)+\frac{1}{\varepsilon}\left(\proj_{\calb}\big(\hat{p}\right)-\hat{p}\big),\\
		\dot{k}&=\tau\Big(\sum_{i\in\cali}x_i-x^*(t)\Big)^\top{\dot{x}^*(t)}-\tau\,\sigma_s\big(|k|\big) k. 
	\end{align*} 
	For the design parameters of the mechanism, we set $ \varepsilon=2\times 10^{-5} $, $ \sigma=10^5 $, and $ k_0=10 $. Noting the bounds of $ a_i $'s, i.e.,  $ 0.004\leq a_i\leq 0.006 $,  the chosen parameters belong to the intervals defined in \eqref{e:adap-par}. 
	Fig. \ref{fig:p-gamma-k-tem} presents the simulation results for $ \tau=1 $, $ \hat{p}(0)=p_0+0.06\bone_n\notin \calb $, and $ k(0)=0 $. The results demonstrate that the price prediction enters the 
	ball $ {B}{(p_0,\rho)} $ and  
	the trust variables converge to one.  Subsequently, the price prediction converges to $ p^*(t) $,  and as a consequence, the aggregative behavior converges to  the desired one as depicted in Fig. \ref{fig:x-sum-tem}.

	\section{Conclusions}\label{sec:con}
	We have presented a nudge framework where a regulator can steer the aggregative behavior of a set of price-taking agents to a desired behavior by sending a suitable price prediction signal. Due to the discrepancy between the signal sent out by the regulator and the actual price, we have incorporated trust dynamics in the agents' model, where the trust variables get updated 
	based on the history of the accuracy of the price prediction signal. Nudge mechanisms have been proposed to steer the aggregative behavior of the agents to desired stationary as well as temporal behaviors. Analytical convergence guarantees have been provided for the proposed nudge mechanisms and the results are demonstrated on a numerical case study.  Future works include investigating the application of the proposed nudge framework in transportation as well as power networks. 
	\bibliographystyle{plain}        
	\bibliography{MyReferences}
	
	\appendix
	
	\textbf{Technical Appendices}
	\section{Existence of solutions for nonautonomous projected dynamical systems}\label{app:pds-exist} 
	\begin{lemma}\label{lem:exis-pds}
		Consider a nonempty compact set $ \calx\subset \R^n $ and a vector field $ h:\R^n\times[0,\infty)\to \R^n$ that is locally Lipschitz  
		in the first argument and measurable in the second. Then, for any initial condition $ x(0)\in\calx $, there exists a Carath\'{e}odory solution $ t\mapsto x(t) $  of the nonautonomous projected dynamical system 
		\begin{equation}\label{non-pds}
			\dot{x}=\Pi_{\calx}\left(x,h(x,t)\right)
		\end{equation}
		satisfying $ x(t)\in\calx $ for all $ t\in[0,\infty) $.
	\end{lemma}
	\BP
	The proof involves 
	demonstrating the existence of  \textit{Krasovskii solutions} for \eqref{non-pds} and then establishing the equivalence of the set of Krasovskii and Carath\'{e}odory solutions.  Since  $ \calx $ is a compact set, we have  the function $ (x,t) \mapsto h(x,t) $ is Lipschitz  on the set $\calx$    \cite[Ex. 3.19]{khalil2002nonlinear} and measurable in $t$. Consequently, by  \cite[Thm. 2]{hauswirth2018time}, the  system admits {Krasovskii solutions}.  Note that in the referred results, the map $h$ is required to be Lipschitz  everywhere in the domain. However, the implication holds even when $h$ is Lipschitz  only on the set $\calx$, that is, the set where the solutions are restricted to. The proof concludes by using \cite[Thm. 6.3]{hauswirth2018projected} which shows that  the set of Krasovskii and Carath\'{e}odory solutions are equivalent for autonomous projected dynamical system. The result extends to the nonautonomous case using the same reasoning. 
	\EP 
	
	\begin{lemma}\label{lem:exis2}
		Consider a nonempty compact set $ \caly\subset \R^m $ and two vector fields $ h:\R^n\times \R^m\times [0,\infty)\to \R^n$ and $ g:\R^n\times \R^m\times [0,\infty)\to \R^m $ that are locally Lipschitz 
		in the first two arguments and measurable in the third one. Consider the nonautonomous projected dynamical system 
		\begin{equation}\label{non-pds2}
			\begin{aligned}
				\dot{x}&=h(x,y,t),\\
				\dot{y}&=\Pi_{\caly}\left(y,g(x,y,t)\right).
			\end{aligned}
		\end{equation} 
		Moreover, assume that there exist a continuously differentiable function $ V:\R^n \to \R $ satisfying:
		\begin{enumerate}[(i)]
			\item $  V(x)\geq 0 $ for all $ x\in \R^n $,
			\item $V(x)\to \infty$ as $\|x\|\to \infty  $,
			\item there exists a constant $ \mu>0 $ such that the following  holds for all $ y\in \caly $, $ t\in[0,\infty) $, and $\|x\|\geq \mu$,
			\begin{equation*}
				\nabla V(x)^\top h(x,y,t)\leq 0. 
			\end{equation*}
		\end{enumerate} 
		\noindent Then, for any initial condition $ (x(0),y(0))\in \R^n\times\caly $, there exists a bounded Carath\'{e}odory solution $ t\mapsto (x(t),y(t)) $  of the system \eqref{non-pds2} over the domain $ [0,\infty) $.
	\end{lemma}
	\BP 
	Our proof proceeds in two steps. First, for each initial condition, we design a nonautonomous projected dynamical system that admits a solution starting from the said initial point. Second, we show that this solution is also a solution of~\eqref{non-pds2}. 
	
	Consider the  continuous and  increasing function $ \alpha(s):=\sup_{\|x\|\leq s} V(x) $ for $ s\geq 0 $. Then, from condition (i) imposed on $V$, we have  
	\begin{equation}\label{lyap:alpha}
		0\leq V(x)\leq \alpha(\|x\|),\quad \forall x\in\R^n.
	\end{equation}		
	Let $(x_0,y_0) \in \R^n \times \caly$ be any initial condition. Define $ \calx_0:=\{x\in\R^n \mid V(x)\leq c\} $ where $ c>\max\{V(x_0),\alpha(\mu)\} $. Then $ x_0\in \inter(\calx_0) $ and  
	the closed ball $\bar{B}(\bze_n,\mu)  $ is in the interior of $ \calx_0 $ as a consequence of \eqref{lyap:alpha}. The former fact follows from $ V(x_0)< c $, and we show the latter  by contradiction. Assume that $\bar{B}(\bze_n,\mu)$ is not in the interior of $ \calx_0 $, then there exists some point $ z_0\in  \bar{B}(\bze_n,\mu)$ such that $ V(z_0)=c $. Since $ \alpha(\cdot) $ is an increasing function, it follows from $ z_0\in  \bar{B}(\bze_n,\mu)$ that  $ \alpha(\|z_0\|)\leq \alpha(\mu) $. Bearing this and $ V(z_0)=c>\alpha(\mu) $ in mind, we have $ V(z_0)>\alpha(\|z_0\|) $ which is in contradiction to   \eqref{lyap:alpha}.  
	Note that (ii) implies that  $ \calx_0 $ is compact. Having defined this set, we now consider a compact set $  \calx $ such that $ \calx_0\subset \inter(\calx) $  
	and introduce the following projected dynamical system 
	\begin{equation}\label{full-psd}
		\begin{aligned}
			\dot{x}&=\Pi_{\calx}\left(x,h(x,y,t)\right),\\
			\dot{y}&=\Pi_{\caly}\left(y,g(x,y,t)\right).
		\end{aligned}
	\end{equation}
	From Lemma~\ref{lem:exis-pds}, this system admits a bounded Carath\'{e}odory solution  $t\mapsto (\hat{x}(t),\hat{y}(t)) $ over the domain $ [0,\infty) $ starting from the chosen initial condition $(x_0,y_0)$. That is, here $(\hat{x}(0),\hat{y}(0)) = (x_0,y_0)$.  
	We next show that this solution $(\hat{x}(\cdot),\hat{y}(\cdot))$ is   
	also a solution of the system \eqref{non-pds2}. Since $ (x_0,y_0)\in \R^n\times \caly $ is chosen arbitrary, this  concludes the proof.
	
	Noting that $ x_0\in \inter(\calx_0) $, the solution $ \hat{x}(\cdot) $ is continuous, and $ \calx_0 $ is compact, there exists some finite time $ T>0 $ such that $ \hat{x}(t)\in \calx_0 $ for all $ t\in [0,T] $. In this time interval, the projection in the $x$-component of~\eqref{full-psd} is not active since $ \calx_0 \subset \inter( \calx) $, that is, we have $  \Pi_{\calx}\left(x,h(x,y,t)\right)= h(x,y,t)$. Bearing this in mind together with  (iii) and $ \bar{B}(\bze_n,\mu) \subset \inter(\calx_0) $, we deduce  that $ \hat{x}(t)\in \calx_0 $ for all $ t\in [0,\infty) $ since $ \dot{V}(x)\leq 0 $ on the boundary of $ \calx_0 $. This implies that the projection operator $ \Pi_{\calx}(x,\cdot) $ is inactive for all times because $ \hat{x}(t) $ is in the interior of $ \calx $. Thus, we conclude that  $t\mapsto (\hat{x}(t),\hat{y}(t)) $ is also a solution of the system \eqref{non-pds2}.
	\EP
	
	\section{Existence of solutions for adaptive nudge}\label{apndx:uub} 
	\begin{lemma}\label{lem:adap-uub}
		Consider the closed-loop system formed by agents' model \eqref{agent} and the adaptive nudge mechanism  \eqref{adap_nudge} with   $t\mapsto x^*(t)$ satisfying Assumption~\ref{asm:xs-t}. Let the design parameters satisfy $ \sigma\in \cali_{\sigma} $, $ k_0\in\cali_{k_0} $, and $ \varepsilon \in \cali_{\varepsilon}$ with the intervals $ \cali_{\sigma} $, $ \cali_{k_0} $, and $ \cali_{\varepsilon} $ given by \eqref{e:adap-par}. Then, for any initial condition $ (\hat{p}{(0)},K(0),\col( \gamma_i{(0)}))\in \R^n\times \R^{n\times n} \times [0,1]^{N}$, there exists a bounded Carath\'{e}odory solution  $t \mapsto \left(\hat{p}(t),K(t),\col (\gamma_i(t))\right)$ of the closed-loop system over the domain $ [0,\infty) $. Moreover, there exist some constants $ \bar p>0 $, $ \bar k>0 $, and a finite time $ T\geq 0 $ such that   we have $ \|\hat{p}(t)\|\leq \bar{p}$ and $ \|K(t)\|_{\operatorname{F}} \leq \bar{k} $ for all $ t\in[T,\infty) $. 
	\end{lemma}
	\BP
	The proof is divided in two parts. The first part focuses on establishing existence of Carath\'{e}odory solutions and the second part shows their ultimate boundedness.  
	
	\textit{Existence of solutions:} 	We use the expression of $ x_i $ given by \eqref{opt_action_cont} and \eqref{pi_def} to rewrite the  adaptive nudge mechanism  \eqref{adap_nudge} as follows: 
	\begin{equation}\label{adap-er}
		\begin{aligned}
			\dot{\hat{p}}&=-\Big(\frac{1}{\varepsilon}I_n+\sum_{i\in\cali}\gamma_i  Q_i^{-1}\Big)\,d(\hat{p})+K\dot{x}^*(t)+\nu(t),\\
			\dot{K}&=\tau\Big(\!-\sum_{i\in\cali}\gamma_iQ_i^{-1}d(\hat{p})+\nu(t)\Big)\dot{x}^{*\top}\!(t)-\tau\sigma_s\big(\|K\|_\text{F}\big) K,
		\end{aligned}
	\end{equation}
	where 
	$ d(\hat{p}):=\hat{p}-\proj_{\calb}(\hat{p}) $ and $\nu(t):=\sum_{i\in\cali}(c_i+\gamma_{i}Q_i^{-1}(\hat\lambda_i-\proj_{\calb}(\hat{p})))-\sum_{i\in\cali}Q_i^{-1}\hat\lambda_i-x^*(t).$
	Note that the term $ \nu(t) $ is bounded for all $ \hat{p}\in \R^n $, $ \gamma_i\in [0,1] $, and $ t\geq 0 $. More precisely,	using  $ \proj_{\calb}\left(\hat{p}\right)\in  \calb$ and boundedness of $ x^*(t) $, there exist some finite $ \bar \nu>0 $ such that we have $ \|\nu(t)\|\leq \bar\nu  $ for all $(\hat{p},\col(\gamma_{i}))\in \R^n\times [0,1]^N$ and $t\geq 0$.
	
	Next, we rewrite the dynamics of the overall closed-loop system in a suitable form to argue existence of solutions. Let $ \varphi:=\vect(K) $	and $ \xi:=\col(\hat{p},\varphi) $, then the closed-loop system, made of \eqref{agent} and \eqref{adap-er}, becomes
	\begin{equation}\label{adap_cl_compact}
		\begin{aligned}
			\dot{\xi}&=h(\xi,\col(\gamma_{i}),t),\\
			\dot\gamma_i&=\Pi_{[0,1]}\left(\gamma_i,\eta_i \psi_i(\|p(t)-\hat{p}\|)\right),\quad \forall i\in\cali,
		\end{aligned}
	\end{equation}
	where $ h $ defines 
	the right-hand side of \eqref{adap-er}.  
	Note that the map $ t \mapsto h(\xi,\col(\gamma_{i}),t)$ is measurable as a consequence of Assumption \ref{asm:xs-t}.  
	Further, using the fact that $\sigma_s$ is Lipschitz and following arguments analogous to those provided in the  
	proof of Theorem \ref{thm:soft}, 
	we deduce that  the map $ (\xi,\col(\gamma_{i}),t) \mapsto h(\xi,\col(\gamma_{i}),t)$ is locally Lipschitz  in $ (\xi,\col(\gamma_{i})) $. 
	Also,  the map $ (\hat{p},t) \mapsto \psi_i(\|p(t)-\hat{p}\|) $ is  locally Lipschitz in $\hat{p}$ and measurable in $t$.  Hence, the  existence of bounded solutions  
	over the domain $ [0,\infty) $ follows from verifying that the hypotheses (i)-(iii) of Lemma~\ref{lem:exis2} hold.    
	The rest of the proof achieves this.  
	
	Consider the following Lyapunov  candidate 
	$
	V(\xi):=\frac{1}{2}\|d(\hat{p})\|^2+\frac{1}{2\tau}\|\varphi\|^2.
	$
	Analogous to the proof of Theorem \ref{thm:soft}, we deduce from Danskin's Theorem that $ \|d(\hat{p})\|^2 $ is differentiable and $ \nabla \|d(\hat{p})\|^2=2\,d(\hat{p}) $. Thus, the function $ V $ satisfies the hypotheses (i) and (ii) of Lemma \ref{lem:exis2}. Our next step is to analyze the inner product of $ \nabla V $ and the function $ h $  given by \eqref{adap_cl_compact}. Hence we define 
	\begin{equation*}
		H(\xi,\col(\gamma_{i}),t):=\nabla V(\xi)^\top h(\xi,\col(\gamma_{i}),t).
	\end{equation*}
	In the following discussion, we show existence of some $ \mu>0 $ such that
	\begin{equation}\label{H-ineq}
		H(\xi,\col(\gamma_{i}),t)\leq 0,\quad \forall \|\xi\|\geq \mu,
	\end{equation}
	for all $ \col(\gamma_{i}) \in [0,1]^N $ and $ t\geq 0 $. This verifies that Lemma \ref{lem:exis2}(iii) holds and establishes existence.  
	
	For simplicity of presentation, we compute $ H $ in the coordinates of $ (\hat{p},K,\col(\gamma_{i})) $. Note that in this coordinates, the Lyapunov candidate becomes $ V(\hat{p},K)=\frac{1}{2}\|d(\hat{p})\|^2+\frac{1}{2\tau}\|K\|_\text{F}^2 $. This allows us to find the relation of  $ H $ as follows 
	\begin{equation*}
		H(\hat{p},K,\col(\gamma_{i}),t)=\Tr\left(\begin{bmatrix}
			\dot{\hat{p}} & \dot{K}
		\end{bmatrix}\begin{bmatrix}
			d(\hat{p})^\top \\ \frac{1}{\tau}K^\top
		\end{bmatrix}\right),
	\end{equation*}
	where $ \begin{bmatrix}
		\dot{\hat{p}} & \dot{K}
	\end{bmatrix} $ stands for the right-hand side of \eqref{adap-er}. Expanding on the expression, we get  
	\begin{equation}\label{vdot-1}
		\begin{split}
			&H=-\|d(\hat{p})\|_{\sum_{i\in\cali}\gamma_iQ_i^{-1}}^2-\frac{1}{\varepsilon}\|d(\hat{p})\|^2\\
			&+d(\hat{p})^\top K\dot{x}^*(t)+d(\hat{p})^\top \nu(t)+\frac{1}{\tau}\Tr\left(\dot K  K^\top\right),
		\end{split}
	\end{equation}
	where
	\begin{equation*}
		\begin{aligned}
			\frac{1}{\tau}\Tr\left(\dot K  K^\top\right) 
			=&-d(\hat{p})^\top\sum_{i\in\cali}\gamma_i  Q_i^{-1} K \dot{x}^{*}(t)\\
			&+\nu(t)^\top K \dot{x}^{*}(t)-\sigma_s\big(\|K\|_\text{F}\big)\|K\|_\text{F}^2.
		\end{aligned}
	\end{equation*}
	In~\eqref{vdot-1}, we have dropped the arguments of $H$ for simplicity.
	Since $ \gamma_{i}\in [0,1] $ and $ Q_i\succ 0 $ for all $ i\in \cali $, the first term on the right-hand side of \eqref{vdot-1} is nonpositive. Hence, we have
	\begin{align*}
		H\leq &-\frac{1}{\varepsilon}\|d(\hat{p})\|^2+d(\hat{p})^\top\Big(I_n-\sum_{i\in\cali}\gamma_i  Q_i^{-1}\Big)K\dot{x}^*(t)\\
		&+d(\hat{p})^\top \nu(t)+\nu(t)^\top K \dot{x}^{*}(t)-\sigma_s\big(\|K\|_\text{F}\big)\|K\|_\text{F}^2\,.
	\end{align*} 
	Further, one can   
	show that $ \|I_n-\sum_{i\in\cali}\gamma_i  Q_i^{-1}\|\leq 1+\lambda_{\max}(\sum_{i\in\cali}Q_i^{-1}) $. This yields   
	\begin{multline*}
		d(\hat{p})^\top\Big(I_n-\sum_{i\in\cali}\gamma_i  Q_i^{-1}\Big)K\dot{x}^*(t) 
		\\
		\leq \frac{\theta}{2}\Big(1+\lambda_{\max}(\sum_{i\in\cali}Q_i^{-1})\Big) \left(\|d(\hat{p})\|^2+\|K\|_\text{F}^2\right),
	\end{multline*}
	where we used $ \|\dot{x}^*(t)\|\leq \theta $ (cf. Assumption \ref{asm:xs-t}),  $ \|K\|\leq \|K\|_\text{F} $ and Young’s inequality  $ 2\|d(\hat{p})\|\|K\|_\text{F}\leq \|d(\hat{p})\|^2+\|K\|_\text{F}^2 $.  Consequently, using the above inequality and the  bounds on $\nu(t) $ and $ \dot{x}^*(t) $, we deduce that
	\begin{equation}
		\begin{aligned}
			H\leq& -\frac{1}{\varepsilon}\|d(\hat{p})\|^2+\frac{\theta}{2}\Big(1+\lambda_{\max}(\sum_{i\in\cali}Q_i^{-1})\Big)\\
			&\cdot \left(\|d(\hat{p})\|^2+\|K\|_\text{F}^2\right)+\bar\nu\|d(\hat{p})\|+\bar\nu \theta\|K\|_\text{F} \\ 
			&-\sigma_s\big(\|K\|_\text{F}\big)\|K\|_\text{F}^2. \label{eq:v-ineq}
		\end{aligned}
	\end{equation} 
	We proceed the proof by showing that, by selecting the design parameters carefully, there exists a compact set such that the right-hand side of the foregoing equation is negative  
	outside of this set. 
	Toward this end,   
	we make use of the definition of $ \sigma_s(\,\cdot\,) $ and deduce that, for any $ \sigma>0 $ and $ k_0>0 $, the last term on the right-hand side of~\eqref{eq:v-ineq} satisfies 
	$
	-\sigma_s\big(\|K\|_\text{F}\big)\|K\|_\text{F}^2\leq -\frac{\sigma}{2}\|K\|_\text{F}^2+\frac{\sigma}{2}k_0^2\,.
	$
	This implies that  
	\begin{align*}
		H\leq& -\frac{1}{\varepsilon}\|d(\hat{p})\|^2-\frac{\sigma}{2}\|K\|_\text{F}^2+\frac{\theta}{2}\Big(1+\lambda_{\max}(\sum_{i\in\cali}Q_i^{-1})\Big)\\
		&\cdot \left(\|d(\hat{p})\|^2+\|K\|_\text{F}^2\right)
		+\bar\nu\|d(\hat{p})\|+\bar\nu \theta\|K\|_\text{F} +\frac{\sigma}{2}k_0^2.
	\end{align*}
	Let $ \varepsilon\in \cali_{\varepsilon} $ and $ \sigma\in\cali_{\sigma} $ with $ \cali_{\varepsilon} $ and $ \cali_{\sigma} $ given by \eqref{e:adap-par}. Then we get
	\begin{equation}\label{H-final}
		\begin{aligned}
			H&\leq
			-\frac{1}{2\varepsilon}\|d(\hat{p})\|^2-\frac{\sigma}{4}\|K\|_\text{F}^2+\bar\nu\|d(\hat{p})\|+\bar\nu \theta\|K\|_\text{F} +\frac{\sigma}{2}k_0^2\\ 
			&= -\frac{1}{4\varepsilon}\|d(\hat{p})\|^2-\frac{\sigma}{8}\|K\|_\text{F}^2-\frac{1}{4\varepsilon}(\|d(\hat{p})\|-2\varepsilon \bar\nu )^2\\
			&\quad-\frac{\sigma}{8}(\|K\|_\text{F}-\frac{4}{\sigma}\bar\nu \theta)^2+c,
		\end{aligned}
	\end{equation}
	where $ c:=\frac{2}{\sigma}\bar\nu^2 \theta^2+\varepsilon \bar{\nu}^2+\frac{\sigma}{2}k_0^2  $. Note that the third and forth terms on the right-hand side of the equality are nonpositive. Consequently, bearing the definition of  $ \xi $ in mind, we obtain \eqref{H-ineq} with $ \mu=\max\{\bar{\delta}+\|p_0\|+2\sqrt{\varepsilon c}\,,\sqrt{8\sigma^{-1}c}\} $. Thus, existence  of the solutions for all $ t\geq 0 $ is guaranteed.

	\textit{Deriving ultimate bounds:} 	Noting $ \varepsilon\in \cali_{\varepsilon} $ and $ \sigma\in\cali_{\sigma} $, we deduce  from \eqref{H-final} that the time-derivative of the evolution of $V$ along any solution of \eqref{adap-er} satisfies
	$
	\dot{V}\leq -\beta V+b
	$ with 
	\begin{align*}
		\beta&= \frac{\theta(1+\lambda_{\max}(\sum_{i\in\cali}Q_i^{-1}))}{2\max\{1,\tau^{-1}\}},\\
		b&=\frac{2}{\sigma}\bar\nu^2 \theta^2+\bar{\nu}^2\theta^{-1}(1+\lambda_{\max}(\sum_{i\in\cali}Q_i^{-1}))^{-1} +\frac{\sigma}{2}k_0^2\,.
	\end{align*} 
	This implies that $ \dot{V}\leq -\frac{\beta}{2} V $ whenever $ V\geq \frac{2b}{\beta} $.
	Thus, along the solution, we have $ V(t)\leq \exp(-\frac{\beta}{2} t) V(0) $ whenever $ V(t)\geq \frac{2b}{\beta} $. It follows that for a solution starting outside of the compact set $ \Omega:=\{(\hat{p},K)\in \R^n\times \R^{n\times n}\mid V(\hat{p},K)\leq \frac{2b}{\beta}\}  $, it converges exponentially fast to $ \Omega $ in the time interval $ [0,T] $ with $T= \frac{2}{\beta}\ln(\frac{\beta V(0)}{2b})$, and remains there afterwards. In addition, for a solution starting in  $ \Omega $, the inequality $ V(t)\leq  \frac{2b}{\beta} $ is satisfied for all $ t\geq  T=0 $ since $ \dot{V} $ is negative on $ \bd(\Omega) $. We  conclude from this argument  that  $ (\hat{p}(t) , K(t)) $ belongs to the set $ \{(\hat{p},K)\in \R^n\times \R^{n\times n}\mid \|\hat{p}\|\leq \bar{p}, \|K\|_{\operatorname{F}} \leq \bar{k}\} $ for all $ t\geq T $, where
	\begin{equation}\label{uub}
		\begin{aligned}
			\bar p &:=\|p_0\|+\bar{\delta}+2\sqrt{\beta^{-1}b},\\
			\bar k &:=2\sqrt{\tau\beta^{-1}b}\,.
		\end{aligned}
	\end{equation}
	\EP
	
\end{document}